%% file: ff.tex
\documentclass[10pt,a4paper,oneside]{article}

\usepackage{latexsym}
\usepackage{amsbsy}
\usepackage{amsmath}
\usepackage{amssymb}
\usepackage{graphics}
\usepackage{epsfig}

\def\RR{\mathbb{R}}
\def\ZZ{\mathbb{Z}}

\def\i{\mathrm{int}}

\newtheorem{thm}{Theorem}[section]
\newtheorem{lem}[thm]{Lemma}
\newtheorem{prop}[thm]{Proposition}
\newtheorem{cor}[thm]{Corollary}

\newenvironment{defin}
      {\vspace{2.5mm}\par\noindent
                \textbf{Definition.}}
             {\vspace{2.5mm}}

\newenvironment{proof}%
      {\par\noindent%
            \textbf{Proof.}}%
           { ~\hfill$\Box$\linebreak}

\newcommand{\dv}{ \! : \! }

\title{\large \bfseries TRIANGULATIONS OF FIBRE-FREE HAKEN 3-MANIFOLDS}

\author{\normalsize ALEKSANDAR MIJATOVI\'{C}}

\date{}

\renewenvironment{abstract}
            {\begin{quotation}\noindent\small
                              \textsc{Abstract}.\hspace{0.5mm}}
            {\end{quotation}}

\begin{document}

\maketitle

\begin{abstract}
It is not known whether there exists a computable function 
bounding
the number of Pachner moves needed to connect any two triangulation of
a compact 3-manifold. 
In this paper we find an explicit bound of this kind for all
Haken 3-manifolds which contain no fibred submanifolds as 
strongly simple pieces of their JSJ-decomposition.
The explicit formula for the bound is in terms of the number of 
tetrahedra in the two triangulations. This implies a conceptually
trivial algorithm for recognising any non-fibred knot complement
among all 3-manifolds.
\end{abstract}

\begin{center}
\section{\normalsize \scshape INTRODUCTION}
\label{sec:intro}
\end{center}

\input{intro}

\begin{center}
\end{center}
\vspace{-8mm}
\begin{center}
\section{\normalsize \scshape CANONICAL DECOMPOSITIONS OF 3-MANIFOLDS}
\label{sec:jsj}
\end{center}

\input{jsj}

\begin{center}
\end{center}
\vspace{-7mm}
\begin{center}
\section{\normalsize \scshape STATEMENT OF THE MAIN THEOREM}
\label{sec:statement}
\end{center}

\input{statement}

\begin{center}
\end{center}
\vspace{-7mm}
\begin{center}
\section{\normalsize \scshape THE CANONICAL HIERARCHY}
\label{sec:can}
\end{center}

\input{canonical}

\begin{center}
\end{center}
\vspace{-7mm}
\begin{center}
\section{\normalsize \scshape PROOF OF THE MAIN THEOREM}
\label{sec:proof}
\end{center}

\input{proof}

\bibliographystyle{amsplain}
\bibliography{cite}

\noindent \textsc{Department of Pure Mathematics and Mathematical Statistics},
\textsc{Center for Mathematical Sciences},
\textsc{University of Cambridge}\\
\textsc{Wilberforce Road},
\textsc{Cambridge, CB3 0WB},
\textsc{UK}\\
\textit{E-mail address:} \texttt{a.mijatovic@dpmms.cam.ac.uk}

\end{document}

%% file: intro.tex
It is a non-trivial fact (proved by Pachner in~\cite{pachner}) that
any triangulation of a
compact PL
$n$-manifold
can be transformed into any other triangulation of the same manifold
by a finite sequence of simplicial moves
and simplicial isomorphisms. The moves can be described
as follows.

\begin{figure}[!hbt]
  \begin{center}
    \epsfig{file=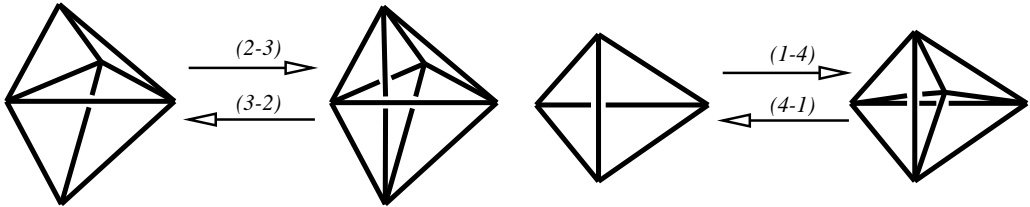}
        \caption{{\small Three dimensional Pachner moves.}}
        \label{fig:3pm}
   \end{center}
\end{figure}

\begin{defin}
Let
$T$
be a triangulation of a compact PL
\mbox{$n$-manifold}
$M$.
Suppose
$D$
is a combinatorial
\mbox{$n$-disc}
which is a subcomplex both of
$T$
and of the boundary of a standard
\mbox{$(n+1)$-simplex}
$\Delta^{n+1}.$
A
\textit{Pachner move}
consists of changing
$T$
by removing the subcomplex
$D$
and inserting
$\partial\Delta^{n+1}-\i(D)$
(for
$n$
equals 3, see figure~\ref{fig:3pm}).
\end{defin}

It is 
an immediate consequence of the definition that there are precisely
$(n+1)$
possible Pachner moves in dimension
$n$.
If our 
$n$-manifold 
$M$
has non-empty boundary, then the moves 
from this definition do not alter the induced triangulation of
$\partial M$.
But changing the simplicial structure 
(throughout this paper the term
\textit{simplicial structure} 
will be used as a synonym for a triangulation) 
of the boundary with an
$(n-1)$-dimensional Pachner move can be achieved by
gluing onto (or removing from) our manifold
$M$
the standard
$n$-simplex
$\Delta^n$
that exists by the definition of the move.
Our aim in this paper is to consider the triangulations of Haken
3-manifolds. 
Their JSJ-decompositions (see section~\ref{sec:jsj} for the
precise definition) consist
of strongly simple pieces,
$I$-bundles
and
Seifert fibred spaces. The strongly simple pieces are the ones
that contain all the interesting topological information
about the manifold
and also have the crucial property
of being atoroidal (i.e.
all incompressible
tori in them are boundary parallel) and are 
hence hyperbolic. It is precisely the strongly simple 
submanifolds that we have to make additional hypothesis on in the next
theorem.\\ 

\begin{thm}
\label{thm:shit}
Let
$M$
be a Haken 3-manifold that does not
contain strongly simple pieces
which are surface bundles or semi-bundles in its JSJ-decomposition.
Let
$P$
and
$Q$
be two triangulations of
$M$
that contain
$p$
and
$q$
tetrahedra respectively. Then there exists a sequence of
Pachner moves of length at most
$e^{2^{ap}}(p)+ e^{2^{aq}}(q)$
which transforms
$P$
into a triangulation isomorphic to
$Q$.
The constant
$a$
is bounded above by
$200$.
The homeomorphism of
$M$,
that realizes this simplicial isomorphism, is supported in
the characteristic submanifold
of
$M$
and it does not permute the components of
$\partial M$.
\end{thm}

The triangulations appearing in theorem~\ref{thm:shit} are 
allowed to be non-combinatorial, which means that the simplices
are not (necessarily) uniquely determined by their vertices.       
The exponent in the above expression containing the exponential
function
$e(x)=2^x$
stands for the composition of the function with itself rather 
than for multiplication. Since the formula in theorem~\ref{thm:shit}
is explicit, it gives a conceptually trivial algorithm (see 
proposition 1.3 in~\cite{mijatov}) to recognise
any 3-manifold that satisfies the hypothesis of the 
theorem (just make all possible
sequences of Pachner moves whose length is smaller than the bound!). 

The unthinkable magnitude of this bound should, I suppose,
be measured against the vastness of the class of 3-manifolds
it covers. It for example gives a direct way of determining
whether any 3-manifold is homeomorphic to a given non-fibred
knot complement. 
In section~\ref{sec:statement} we outline a simple
procedure, based on theorem~\ref{thm:non-fibred}
(which is a slight generalisation of~\ref{thm:shit}), that
can be used to decide if a knot, represented by a knot
diagram, is the same as our given non-fibred knot.

\vspace{-2mm}
The proof of theorem~\ref{thm:non-fibred} is carried
out in several phases and uses a variety 
of techniques. We start by subdividing the original triangulation
so that the characteristic submanifold is supported by a subcomplex
of the subdivision. Then we apply the canonical hierarchy techniques
(see section~\ref{sec:can}) 
together with theorem 1.2 of~\cite{mijatov} to 
the strongly simple pieces. 
This makes it 
possible to connect any two triangulations of the simple submanifolds
by a sequence of Pachner moves. 
The subdivision of the original triangulation in the
characteristic submanifold can be altered directly by applying
the main theorem 
of~\cite{mijatov1}.

In section~\ref{sec:jsj} we give a brief exposition of JSJ-theory,
generalise it to the setting of 3-manifolds with boundary pattern
and prove a bound on the normal complexity of the canonical surfaces. 
Section~\ref{sec:statement} gives a precise definition of fibre-free
3-manifolds and states our main theorem. In section~\ref{sec:can}
we define the canonical hierarchy and prove that it has all the required 
properties. The last section puts everything together and proves 
theorem~\ref{thm:non-fibred}.

%% file: jsj.tex
Any 3-manifold contains a (possibly empty) collection
of canonical tori and annuli. Cutting along 
these surfaces we obtain
a canonical decomposition of our space.
These, so called, JSJ-decompositions 
of 3-manifolds 
are 
due to Jaco-Shalen~\cite{js} and Johannson~\cite{joh} 
with ideas from Waldhausen.
When studying a triangulation of a Haken 3-manifold,
it is profitable to make it interact well with the pieces of
the JSJ-decomposition. In other words the first step towards
simplifying the triangulation of our 3-manifold will consist of
subdividing the original triangulation so that the pieces of the
JSJ-decomposition are triangulated by the subcomplexes of the subdivision.

In subsection~\ref{subsec:char} we are going to define canonical surfaces, 
JSJ-decompositions and characteristic 
submanifolds (main reference for this subsection is~\cite{neumann}). Then,
in~\ref{subsec:Pcan}, we will
study the parallel theory of canonical annuli in the presence of
boundary patterns. This will be very useful when analysing the 
topology of hyperbolic pieces of our 3-manifold (see section~\ref{sec:can}).
And finally in~\ref{subsec:csnf}
we will construct the canonical tori
and annuli from fundamental surfaces. 

\vspace{-5mm}
\begin{center}
\end{center}
\begin{center}
\subsection{\normalsize 
            \scshape JSJ-DECOMPOSITION AND THE CHARACTERISTIC SUBMANIFOLD}
\label{subsec:char}
\end{center}

Consider an irreducible 3-manifold
$M$
with (possibly empty) incompressible boundary. 
Recall that a properly
embedded annulus or torus in 
$(M,\partial M)$
is called \textit{essential} if it is incompressible
and not boundary parallel. The manifold
$M$
is \textit{atoroidal} (resp. \textit{an-annular}) if it
contains no essential tori (resp. annuli).
An essential annulus or torus 
$S$ 
is called
\textit{canonical}
if any other essential annulus or torus in
$M$
can be isotoped to be disjoint from 
$S$.

Now we want to look at a disjoint maximal collection
$\{S_1,\ldots, S_k\}$
of canonical surfaces in
$M$
such that no two of the
$S_i$
are parallel.
Since we are assuming that 
$\partial M$
is incompressible in 
$M$,
the essential annuli in
$\{S_1,\ldots, S_k\}$
are also boundary incompressible. So a straightforward
application of the Kneser-Haken finiteness theorem (see theorem III.20 
in~\cite{jaco1})
guarantees the existence
of such a maximal collection. The result
of cutting 
$M$
along such a union of canonical surfaces is
sometimes referred to as a \textit{Waldhausen decomposition}
of
$M$.
It is shown in~\cite{neumann} (see lemma 2.2) that
a maximal system of disjoint canonical surfaces 
$\{S_1,\ldots, S_k\}$
is unique up to isotopy.
In other words any incompressible annulus or torus
$S$
in
$(M,\partial M)$
can be isotoped to be disjoint from the
surface
$S_1\cup\ldots\cup S_k$.
Moreover, if
$S$
is not parallel to any 
$S_i$,
then its final position in
$M-\mathrm{int}\mathcal{N}(S_1\cup\ldots\cup S_k)$
is determined up to isotopy.

Let's now look at a 
piece 
$M'$
of the Waldhausen decomposition of
$M$.
Put differently
$M'$
is simply a component 
of the cut-open manifold
$M-\mathrm{int}\mathcal{N}(S_1\cup\ldots\cup S_k)$.
Let 
$\partial_1M'$
be the part of
$\partial M'$
coming from the surface
$S_1\cup\ldots\cup S_k$
and
let
$\partial_0M'$
equal
$M'\cap \partial M$.
Clearly the union of
$\partial_0M'$
and
$\partial_1M'$
equals
$\partial M'$
and the components of
$\partial_1 M'$
are annuli and tori. Also both surfaces
$\partial_0M'$
and
$\partial_1M'$
are incompressible while it is possible for
$\partial M'$
to compress into 
$M'$.
We say that a piece
$M'$
is
\textit{simple}
if any essential annulus or torus in
$(M', \partial_0M')$
is parallel to a component of
$\partial_1 M'$.
If
$M'$
is a simple piece which does not admit an incompressible
annulus which is properly embedded in
$(M', \partial_1M')$,
then we call it
\textit{strongly simple}.
It turns out that all the pieces that are simple but not 
strongly simple are either Seifert fibred or of the
form
$(\mathrm{torus})\times I$
with 
$\partial_1 M'=(\mathrm{torus})\times \partial I$
(the manifold 
$M$
is in this case homeomorphic to a torus bundle over 
a circle with holonomy of trace different from 
$\pm2$). 
The simple Seifert fibred pieces are of course very restricted
as well (see proposition 3.2 and figure 1 in~\cite{neumann}).
It is the topology of strongly simple pieces from our
3-manifold 
$M$
that we will be exploring in section~\ref{sec:can}.
The central result of the JSJ-theory (proposition 3.2 in~\cite{neumann})
says that each complementary 
piece
$M'$
falls into one of the three categories:
\vspace{2mm}

\begin{itemize} 
\item[(a)] $(M', \partial_0 M')$
           is strongly simple.
\item[(b)] $(M', \partial_0 M')$ 
           is an 
           $I$-bundle over a (possibly closed) surface.
\item[(c)] $(M', \partial_0 M')$
           is a Seifert fibred space with a possibly
           non-orientable base space.
\end{itemize} 

These possibilities are almost mutually exclusive. 
Assuming that 
$\partial_1 M'$ is not empty, 
$M'$
can not be both strongly simple and Seifert fibred.
The only 
$I$-bundles that are also strongly simple are the ones
over twice punctured disc and once punctured Moebius band. The
only 
$I$-bundle that is also Seifert fibred is the one over a Moebius band
(see proposition 3.3 in~\cite{neumann} for these non-uniqueness statements).

A canonical annulus in the 3-manifold
$M$
separating two fibred pieces whose fibrations match along it,
or separating a fibred piece from itself so that the 
fibrations on both sides of the annulus still match, is called a 
\textit{matching annulus}. Notice that a piece
$M'$
of the Waldhausen decomposition
containing a matching annulus has to be Seifert fibred, because
an annulus separating two 
$I$-bundles is never canonical. It also not hard to see that a canonical
torus can not separate two pieces with matching fibrations. However
matching annuli do exist (see lemma 3.4 in~\cite{neumann}), but they
carry no topological information.
By deleting all matching annuli from our disjoint maximal collection
of canonical surfaces
$\{S_1,\ldots, S_k\}$
we obtain the 
\textit{JSJ-system}
of canonical surfaces for 
$M$.
The decomposition of
$M$
along this JSJ-system is called the 
\textit{JSJ-decomposition}
of
$M$.
Notice that the pieces of the JSJ-decomposition of
$M$
still fall into the three categories mentioned above.  

Let's define a submanifold 
$\Sigma$
of 
$M$
in the following way.
Let 
$\Sigma$
be a union of the
$I$-bundle pieces and the Seifert fibred pieces of the
JSJ-splitting of
$M$.
Let it also contain a regular neighbourhood of
every canonical annulus or torus that separates two
strongly simple pieces of the JSJ-decomposition, neither
of which is an 
$I$-bundle.
If two pieces of the JSJ-splitting 
of
$M$,
that are already contained in
$\Sigma$,
meet along
a canonical torus or annulus, we remove the interiors of their
regular neighbourhoods from
$\Sigma$.
The submanifold 
$\Sigma$
defined in this way is called the
\textit{characteristic submanifold} of
$M$.
Notice that the tori and annuli whose regular neighbourhoods we've
removed from 
$\Sigma$
when defining it,
are precisely the ones along which two non-matching fibred pieces
meet (e.g. a Seifert fibred piece and an $I$-bundle). It also follows
directly from the definition that every essential annulus or torus
in
$M$
can be isotoped into the characteristic submanifold 
$\Sigma$.

\vspace{-5mm}
\begin{center}
\end{center}
\begin{center}
\subsection{\normalsize \scshape JSJ-THEORY WITH BOUNDARY PATTERN}
\label{subsec:Pcan}
\end{center}

In this subsection and throughout this paper
it will be assumed that the reader is familiar with
section 2 of~\cite{mijatov1}. In particular this amounts to
the basic definitions of incompressibility and
$\partial$-incompressibility of embedded surfaces as well
as some normal surface theory. 
We are now going to define 
a key concept for the construction of the canonical
hierarchy.

A \textit{boundary pattern} 
$P$
in a compact 3-manifold
$M$
is a (possibly empty) collection of disjoint
simple closed curves and trivalent graphs embedded in
$\partial M$
such that the surface
$\partial M - P$
is incompressible in 
$M$.
Assume from now on that our 3-manifold
$M$
is equipped with a boundary pattern 
$P$.
Boundary patterns appear naturally in the following way.
If
$M$
contains a properly embedded two-sided incompressible surface
$S$
with
$\partial S$
intersecting 
$P$
transversally (and missing the vertices of
$P$),
then the cut-open manifold 
$M_S$
inherits a boundary pattern as follows. If
$S'$
and
$S''$
are the two copies of
$S$
in
$\partial M_S$,
then the new boundary pattern, lying in
$\partial M_S$,
can be defined as
$(P\cap \partial M_S)\cup \partial S' \cup \partial S''.$
Our definition of boundary pattern implies that 
the manifold
$M$
has incompressible boundary if and only if it admits 
an empty boundary pattern. 

We shall use boundary patterns to 
keep track of the topological
information as we move down the canonical hierarchy
(see section~\ref{sec:can}).
This is precisely the idea that
Haken exploited to find an algorithm for classifying non-fibred
3-manifolds that contain an injective surface. 
At the heart of Haken's classification program lies the
concept of a $P$-canonical annulus, that shall be
described shortly.

Let 
$M$
be a 3-manifold with non-empty 
boundary that contains a boundary pattern 
$P$. 
Recall that a subset
of
$M$
is 
called \textit{pure},
if it has empty intersection with the pattern
$P$.
Most concepts from general 3-manifolds carry over to 
3-manifolds with pattern in a very natural way. 
For example 
a properly embedded surface
$F$
in
$M$
is \textit{$P$-boundary incompressible}
if for any pure disc 
$D$
in
$M$,
such that 
$D\cap (\partial M\cup F)=\partial D$
and
$D\cap F$
is a single arc in 
$F$,
the arc 
$D\cap F$
cuts off a pure disc from
$F$.
Notice that our definitions 
imply that a 
pure incompressible annulus 
in
$M$
is $P$-boundary incompressible 
if and only if it
is not parallel to a pure annulus in
$\partial M$.
A $P$-boundary incompressible
pure annulus 
$A$
in
$M$
will be called
\textit{trivial}
if it is parallel 
(rel $\partial A$) 
to an annulus 
in the boundary of
$M$.
The interior of this annulus in
$\partial M$
that 
$A$
is parallel to,
must have a non-empty intersection with the pattern
$P$.
Furthermore 
no spanning arc of 
this annulus in
$\partial M$
can be pure.
All this follows from the definition of the boundary pattern 
$P$
and the fact that
$A$
is $P$-boundary incompressible.
Also, an incompressible 
$P$-boundary incompressible pure annulus in
$M$
is termed 
\textit{$P$-essential}.
So, according to our definitions, a $P$-essential annulus can be trivial.
The next concept is of great importance in all that follows. 

\begin{defin}
A properly embedded annulus 
$A$
in
$M$
is a 
\textit{$P$-canonical annulus} if it 
is  non-trivial and $P$-essential
in
$M$
and if it 
has the following properties:  
\begin{list}{$\bullet$}{\itemsep -1mm\topsep0mm}
\item any 
$P$-essential
annulus in
$M$
can be isotoped off of 
$A$
by an ambient isotopy that is invariant on the pattern
$P$ and 

\item any incompressible torus, that is not parallel to a
(possibly non-pure) boundary component of 
$M$,
can be isotoped off of
$A$
by an isotopy that is fixed on the boundary of
$M$.
\end{list}
\end{defin}

The notion of a $P$-canonical annulus in 3-manifolds
with pattern is a direct generalisation of the notion
of a canonical annulus coming from JSJ-theory. An 
observation that will be useful in section~\ref{sec:can},
is that a solid torus with some pattern on its boundary can 
contain a $P$-essential annulus, but this annulus has to 
be trivial and therefore not $P$-canonical.
Or more generally, a $P$-canonical annulus in 
$M$
can not be parallel to an annulus in the boundary of
$M$.
The next theorem will play a key role in the
construction of the canonical hierarchy (see section~\ref{sec:can}).
Its proof is a good example of how the ideas from 
JSJ-theory generalise naturally 
to 3-manifolds with pattern.
\vspace{1mm}

\begin{thm}
\label{thm:bundle}
Let 
$M$
be an irreducible atoroidal 3-manifold with 
boundary pattern
$P$.
If 
$M$
contains a non-trivial $P$-essential annulus 
but no $P$-canonical one, then either of the 
two possibilities must occur:
\begin{list}{}{\itemsep 1mm\topsep1mm}
\item[(a)]
The manifold
$M$ 
is homeomorphic to an 
$I$-bundle 
over a (possibly non-orientable) compact 
surface
$F$
and its horizontal boundary 
(i.e. 
$\partial I$-bundle over 
$F$)
is pure. The pattern 
$P$
is contained 
in the vertical boundary 
(i.e. in the
$I$-bundle over
$\partial F$)
and no fibre over any point in
$\partial F$
can be pure.
\item[(b)] The manifold 
$M$ 
is an atoroidal Seifert fibred space (for an explicit description
see figure 3 in~\cite{mijatov1}). 
\end{list}
\end{thm}
%
%

\begin{proof}
We begin by taking a maximal collection 
$\{A_1,\ldots, A_k\}$
of disjoint non-parallel $P$-essential annuli in
$M$
(two pure annuli are parallel if they are parallel
in the usual sense and the parallelity region 
between them is pure).
Since we can assume that 
$M$
is triangulated and that the pattern 
$P$
is contained in the 1-skeleton of that triangulation, any 
collection of disjoint $P$-essential annuli can be put into normal form.
So our maximal collection of disjoint $P$-essential
annuli exists by Kneser-Haken finiteness theorem.
Notice also that, by assumption, we can make sure that 
at least one of the annuli
$A_i$
is non-trivial.
Let 
$N$
be a component of the cut-open manifold
$M-\mathrm{int}\mathcal{N}(A_1\cup\ldots\cup A_k).$
We will show that 
$N$
fibres either as an 
$I$-bundle or as a Seifert fibred space. Moreover these
fibrations will match up when we reglue the pieces to form the
manifold 
$M$
we started with. The atoroidallity of
$M$
will impose severe restrictions on the Seifert fibred spaces 
that can arise.
Most of the above statements will follow from the next claim.\\
\noindent \textbf{Claim.} The manifold
$N$
is either homeomorphic to a pure 
$I$-bundle over an annulus, a punctured annulus, a Moebius band,
or a punctured Moebius band or it is homeomorphic to a 
Seifert fibred space over a base surfaces from 
figure~\ref{fig:sf1}. The pattern on the toral component of
$N\cap \partial M$ 
in figure~\ref{fig:sf1}.3
intersects every fibre of the Seifert fibration which is contained
in that torus.
All annular components of
$N\cap \partial M$ 
in all the cases
are pure except when
the exceptional fibre in figure~\ref{fig:sf1}.2 is not singular.
Then the annulus from
$N\cap \partial M$
must contain components of the pattern 
$P$
and no spanning arc of that annulus 
can be disjoint from the pattern. 

\begin{figure}[!hbt]
 \begin{center}
    \epsfig{file=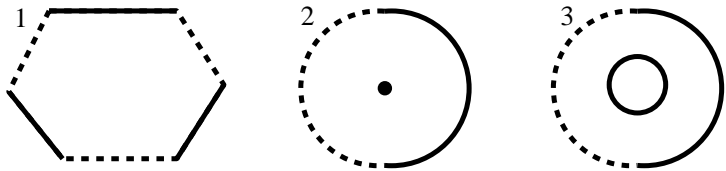}
  \caption{\small The Seifert fibred pieces
                  that can arise as the manifold $N$. 
                  The solid (respectively dashed) part of the boundary
                  of the base surfaces in the figure correspond to the
                  part of the boundary of $N$ coming from 
                  $N\cap \partial M$ (respectively 
                  $N\cap \partial\mathcal{N}(A_1\cup \ldots \cup A_k$)).
                  The dot in case 2 represents a fibre
                  that may or may not be singular.}
                  
  \label{fig:sf1}
 \end{center}
\end{figure}

\vspace{-4mm}
Before proving the claim, we should note that
if
$N$
is a pure 
$I$-bundle over an annulus, then
$M$
has to be an 
$I$-bundle over a torus or a Klein bottle
with empty pattern
and the theorem follows.
This is because 
$M$
contains no $P$-canonical annuli and also no two annuli in
$\{A_1,\ldots, A_k\}$
are parallel.

We are assuming that none of the $P$-essential annuli 
$A_i$
are $P$-canonical. So, an annulus 
$A_i$
is either trivial or there exists another $P$-essential
annulus or an essential torus
$A_i'$
in
$M$,
such that the intersection 
$A_i\cap A_i'$
is minimal and not empty.
The first possibility gives a solid torus as in 
figure~\ref{fig:sf1}.2 with the boundary
as described in the claim. The dot in this case represents a fibre
which is non-singular.

We can assume now that 
$N$
is not a solid torus as in figure~\ref{fig:sf1}.2
with no exceptional fibres (such a piece 
must exist in
$M$
because 
$M$
contains at least one non-trivial
$P$-essential annulus). This means that it contains in its
boundary a copy of an annulus
$A_i$
which is not trivial.
We shall now analyse the 
intersection
$A_i\cap A_i'$
for
the annulus
$A_i$
that is contained in
$\partial N$.
The 1-manifold
$A_i\cap A_i'$
contains no trivial simple closed curves and no
separating properly embedded arcs in either 
of the two surfaces (among other things here,
we are using the fact that the surface
$\partial M - P$
is incompressible in
$M$).
The intersection can therefore consist either of non-trivial 
simple closed curves in
$A_i$
and
$A_i'$,
or of spanning arcs if both of our surfaces
are annuli. The atoroidallity
of
$M$
will impose limitations on the former possibility.

Let's first look at the latter case.
We are therefore assuming that the surface
$A_i'$
is an annulus.
The intersection
$N\cap A_i'$
consists of a disjoint union of discs of length four which inherit a natural
product structure from the annulus
$A_i'$.
In other words they
are of the form
$I\times I$.
Let 
$D\subset N$ 
be one of those discs which intersects
$A_i$,
and let
$s$
be a component of
$D\cap A_i$.
The product structure on
$D$
naturally decomposes the 
$\partial D$
into four arcs,
$s$
being one of them. The two segments
adjacent to 
$s$
are contained in 
$\partial M$
and are disjoint from the pattern 
$P$.
Let 
$t$
be the fourth segment of
$\partial D$,
lying somewhere in the surface
$A_1\cup\ldots\cup A_k$,
say in the component
$A_j$.
We distinguish two subcases depending on whether 
$A_j$
equals
$A_i$
or not.

For the proof of the fact that these 
two subcases yield the 
$I$-bundle possibilities 
in our claim, we refer the reader to
the proof of Case A1 and Case A2 of proposition 3.2
in~\cite{neumann}. It is worth bearing in mind that
if any of the annuli constructed by cut-and-paste 
in the cited proof are compressible, then they have to 
be isotopically trivial (i.e. they bound
a cylinder in
$N$). This is because both discs, obtained by compressing
this
cut-and-paste annulus, chop off two 3-balls from 
$M$.
Neither of these 3-balls can contain 
a properly embedded incompressible annulus.
They can therefore not be nested and our cut-and-paste annulus must be 
isotopically trivial.

Let's now look at the case where
$A_i\cap A_i'$
consists of non-trivial simple closed curves.
The surface 
$A_i'$
can either be a $P$-essential annulus or an essential torus.
This will make 
$N$
into a Seifert fibred space as described by the claim.
The intersection 
$N\cap A_i'$
consists of a disjoint union of annuli. Let
$A$
be one of these annuli that intersects 
$A_i$,
and let 
$s$
be a simple closed curve component of
$A_i\cap A$.
Let 
$t$
be the other boundary component of
$A$.
We have three possibilities according to the position
of the simple closed curve 
$t$.

If 
$t$
is contained in
$N\cap \partial M$,
then cutting 
$A_i$
along
$s$
and pasting parallel copies of 
$A$
to the resulting annuli gives a pair of pure annuli 
that have to be $P$-essential.
They inject into 
$M$
because 
$A_i$
injects. If either of them were 
$P$-boundary compressible, then we could have reduced the number
of components of the intersection
$A_i\cap A_i'$.
They are therefore parallel to components of
$A_1\cup\ldots\cup A_k$
which gives figure~\ref{fig:sf1}.1.

If 
$t$
is contained in
$A_j$
and
$j$
is different from
$i$,
cutting and pasting as before gives two incompressible annuli
that can not both be $P$-boundary compressible (that would make
$A_i$
and
$A_j$
parallel).
They can not both be parallel into
$\{A_1,\ldots , A_k\}$,
because that would contradict the maximality of the family
$\{A_1,\ldots , A_k\}$.
Thus one is parallel into
$N\cap \partial M$
and the other to an annulus from
$N\cap \mathcal{N}(A_1\cup\ldots\cup A_k)$.
This gives figure~\ref{fig:sf1}.1 again.

If 
$t$
lies in
$A_i$,
we have two further subcases coming from
$A$
meeting
$A_i$
from the same side or from the opposite sides.
The latter possibility gives the same as the previous
case, except for one extra case when both annuli produced by
cut-and-paste are
parallel into 
$N\cap \partial M$.
This leads to 
$N$
being an 
$I$-bundle over an annulus. The manifold  
$M$
is then an 
$I$-bundle over a torus or over a Klein bottle with 
empty boundary pattern.

If
$A$
meets 
$A_i$
from the same side and if the orientations on 
$s$
and
$t$
coincide,
then cutting and pasting as before yields
a torus and an annulus. The torus has to be
either boundary parallel or compressible since 
$N$
is atoroidal. The annulus can not be parallel
into
$\{A_1,\ldots, A_k\}$,
because that would give rise to a new $P$-essential annulus
in
$N$,
contradicting the maximality of
$\{A_1,\ldots, A_k\}$.
So the annulus has to be $P$-boundary compressible. This 
gives possibilities 2 and 3 from figure~\ref{fig:sf1}. 
In case of figure~\ref{fig:sf1}.2 the exceptional fibre 
has to be singular, because otherwise we could reduce the 
number of components of the intersection
$A_i\cap A_i'$.
We should also note that the pattern on the toral 
boundary component in figure~\ref{fig:sf1}.3 has to be such
that no fibre of the Seifert fibration on that torus is pure.

The last case, when 
$A$
meets 
$A_i$
from the same side but the orientations
$s$
and
$t$
do not coincide, can not occur. Cut-and-paste as before gives
an annulus which is either parallel into 
$N\cap \partial M$
or
$\{A_1,\ldots, A_k\}$.
In both cases
$N$
is an
$S^1$-bundle over a Moebius band with a part of its
boundary coming from a pure annulus in
$N\cap \partial M$.
We can easily construct a pure annulus in
$N$,
running once around the Moebius band, contradicting the maximality
of
$\{A_1,\ldots, A_k\}$.
This proves the claim.

It is clear that the only manifold pieces
from the claim,
for which the fibration is not unique up to isotopy on their
boundaries, are 
an
$I$-bundle over a Moebius band and an
$I$-bundle over an annulus.
The first one also fibres as a Seifert fibred space from
figure~\ref{fig:sf1}.2 with a singular fibre of index 
$\frac{1}{2}$.
The second one we do not need to consider because we already know 
that the theorem is true if such a piece appears. Since
the non-unique piece
$N$
has a single annulus coming from the collection
$\{A_1,\ldots, A_k\}$,
we can always uniquely extend the fibration on the boundary.
Since we can do that for all other pieces of the complement
$M-\mathrm{int}\mathcal{N}(A_1\cup\ldots\cup A_k)$
as well, it follows that 
$M$
is either an 
$I$-bundle with the described boundary pattern or a Seifert 
fibred space. But since the Seifert fibred space
is also atoroidal, the theorem 
follows.
\end{proof}

\vspace{-10mm}
\begin{center}
\end{center}
\begin{center}
\subsection{\normalsize \scshape CANONICAL SURFACES IN NORMAL FORM}
\label{subsec:csnf}
\end{center}

We are now going to describe how to construct canonical tori and annuli
in a triangulated 3-manifold
using normal surface theory. Our proofs rely on the key notion of a trivial
patch and are based on the following lemma and theorem (see 
section 2 in~\cite{mijatov1}):\\

\begin{lem}
\label{lem:patches}
Let
$M$
be an irreducible 3-manifold with a (possibly empty) boundary
pattern
$P$.
Let
$F$ 
be a minimal weight incompressible P-boundary incompressible
normal surface. If the sum
$F=F_1+F_2$
is in reduced form then each patch is both incompressible and
P-boundary incompressible and no patch is trivial. Furthermore
if
$F$
is injective, then each patch has to be injective.
\end{lem}

\begin{thm}
\label{thm:sum}
Let
$M$
be an irreducible 
3-manifold with a possibly empty boundary pattern
$P$.
Let
$F$ 
be a least weight normal surface properly embedded in
$M$.
Assume also that
$F$
is two-sided incompressible
$P$-boundary incompressible and
$F=F_1+F_2$.
Then
$F_1$
and
$F_2$
are incompressible and
$P$-boundary incompressible.
\end{thm}

Propositions~\ref{prop:cantori} and~\ref{prop:canannuli} are going to
be crucial when subdividing the original triangulation of
$M$.
In particular we will find proposition~\ref{prop:canannuli} 
indispensable while deciphering the topology of the strongly simple pieces
in the JSJ-decomposition of 
$M$.\\

\vspace{-4mm}
\begin{prop}
\label{prop:cantori}
Let 
$M$
be an irreducible 3-manifold with  (possibly empty) 
incompressible boundary. 
Let
$t$
be the number of tetrahedra in the triangulation 
$T$
of
$M$. 
Then every canonical torus in 
$M$
can be isotoped into normal form so that 
it contains not more than
$2^{80t^2}$
normal discs.
The same statement remains true if a single
canonical torus in 
$M$
is replaced by a disjoint maximal collection of 
non-parallel canonical
tori in
$M$.
\end{prop}

\begin{proof}
We start by taking our maximal collection 
$C$
of disjoint canonical 
tori in
$M$
and putting it into normal form, so that
it has the smallest weight in 
its isotopy class. The normal surface 
$C$
can then be expressed as a sum of fundamental surfaces
$C=k_1F_1+\cdots+k_nF_n$
where
$k_1,\ldots,k_n$
are positive integers. 

\noindent \textbf{Claim 1.} Each surface
$F_i$ 
is either an incompressible torus or an incompressible Klein bottle.

We will get the incompressibility of
$F_i$
by applying theorem~\ref{thm:sum}.
All we need to do is to fix a copy of the surface 
$F_i$
and make all regular alterations along curves of
intersection (in the sum
$k_1F_1+\cdots+k_nF_n$)
that do not lie in our copy of 
$F_i$.
This yields a sum
$C=F_i+F'$,
where 
$F'$
is some normal surface in
$M$.
Since both
$C$
and
$M$
are orientable
the surface 
$C$
is two-sided. 
Also, since the normal representative for
$C$
has minimal weight, we can apply theorem~\ref{thm:sum} to conclude
that
$F_i$
must be incompressible as well. 

If we isotope the sum
$C=F_i+F'$
into reduced form, we can apply lemma~\ref{lem:patches} 
to get that no patch can be trivial. 
This means that the surface
$F_i$
can not be a 2-sphere. We are assuming that  
$M$
contains at least one canonical torus and is therefore
different from 
$\RR P^3$.
Since it is also irreducible,
the surface
$F_i$
can not be a projective plane either. Now 
all fundamental surfaces in our sum are closed and connected. 
So the claim follows by the additivity 
of the
Euler characteristic. 

It should be noted that the proof of the claim 
uses only the incompressibility of
$C$
and not the fact that it is canonical. Also the same method of proof
gives that any connected normal surface which appears as a summand of
$C$
has to be an incompressible torus 
or an incompressible Klein bottle.

We now want to bound the coefficients
$k_i$
in the sum
$C=k_1F_1+\cdots+k_nF_n$.
After making all regular alterations in the above sum, except the ones
on the surface
$k_iF_i$,
our expression can be rewritten as
$C=k_iF_i+S$.
The surface 
$S$
is a (possibly disconnected) closed incompressible normal surface with
zero Euler characteristic.

Let's now investigate the patches 
of the polyhedron
$K=F_i\cup S$,
i.e. the surface components of the complement of
the singular locus of
$K$. 
It is clear that our surface 
$C$
lies in a regular neighbourhood of 
$K$.  
The plan is to alter 
$K$
in such a way, so that it contains
no disc patches but it is still a union of 
normal surfaces whose integral
linear combination represents 
$C$.
We want to do this in such a way,
so that the weight of the 
new normal surfaces is not greater than that of
$F_i$
and
$S$. 

If there exists a trivial patch 
$D_i$
in 
$F_i$,
bounded by a simple closed curve from 
$F_i\cap S$,
then this curve must bound a disc 
$D$
in
$S$
which is distinct from
$D_i$. 
We now look at each patch 
$D_i$
of
$K$
that is contained 
in
$F_i$,
such that a copy of the disc
$D$
is also a patch
in
$K$. 
In other words the singular
locus of 
$K$
intersects 
$D$
in its boundary only, making the disc 
$D$
into a patch in the surface
$C$. 

After making all regular alterations in the 
sum
$k_iF_i+S$,
one of the parallel copies of the disc
$D_i$
is adjacent to the disc 
$D$
in the surface
$C$
(see figure~\ref{fig:patch}).
Since the manifold
$M$
is irreducible, the union of discs
$D\cup D_i$
bounds a 3-ball
$B$
as in figure~\ref{fig:patch}. 
If the weights of the discs were not equal, then
we could use this 3-ball 
$B$
to isotope
$C$
so that its weight is reduced. 
Thus we must have the equality
$w(D_i)=w(D)$.

\begin{figure}[!hbt]
 \begin{center}
    \epsfig{file=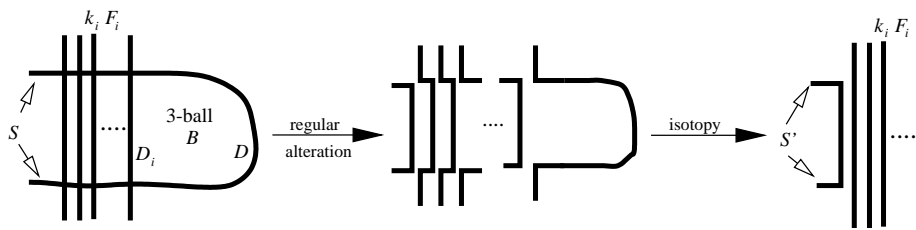}
  \caption{\small Regular alteration along the curve
                  $\partial D_i$ and the isotopy across the 3-ball
                  $B$.} 
  \label{fig:patch}
 \end{center}   
\end{figure}

Now we modify the surface
$C$
by pushing the disc
$D$
across the ball
$B$.
Notice that this operation did not change
the weight of 
$C$.
Also the surface we end up with is again in normal form.
A similar kind of modification can be done to the surface
$S$.
In this case we push the disc
$D$
so that it misses all parallel copies of 
$D_i$
but it is still normally parallel to the disc
$D_i$.
This produces a new normal surface
$S'$
which has the 
same weight as 
$S$.
Also the normal equation
$C=k_iF_i+S'$
holds for the isotoped surface
$C$.

We keep doing the above procedure to all 
patches of
$K$
which are discs
in
$F_i$
and have adjacent trivial patches from
$S$
next to them.
In the end we obtain a polyhedron
$K'=F_i\cup S'$
where
$S'$
is a normal surfaces with the property
$w(S')=w(S)$
and
$K'$
contains no pairs of adjacent disc patches. 

\noindent \textbf{Claim 2.} The polyhedron
$K'$
contains no patches which are discs.

It is clearly enough to show that the sum
$C=k_iF_i+S'$
contains no trivial patches.
Let's assume that there exists a trivial patch 
$R$
somewhere in
$C$
and that it has the smallest weight among all 
trivial patches.
Since
$\chi(F_i)=\chi(S)=0$,
the boundary of
$R$
has to be two-sided in both summands.
So there exists
a unique disc
$R'$
in 
$C$ 
whose boundary is adjacent to that of
$R$.
Disc
$R'$
is different from the trivial patch
$R$,
but it might contain it.
Also
$R'$ 
is not itself a patch, since
$K'$
contains no pairs of adjacent disc patches.
We must have 
$w(R')=w(R)$.
Otherwise we could reduce the weight of
$C$.
$R'$ therefore contains precisely one
trivial patch 
$Q$
with 
$w(Q)=w(R)$
and also other annular patches of zero weight.
Now we take 
the patch
$Q$
and repeat the procedure. We can use the same argument
as in the proof of lemma 3.6 in~\cite{bart}
to express the surface
$C$
as a sum
$C=C'+B$
where 
$B$
is a normal surface with zero weight 
(cf. the proof of lemma~\ref{lem:patches} in~\cite{mijatov1}). 
This contradiction proves the claim.


Claim 2 implies that all patches of
$K'$
are incompressible subsurfaces of
$C$
with non-positive Euler characteristic.
Therefore they have to be 
incompressible annuli. 
So the regular neighbourhood 
${\cal N}(K')$
supports a structure of a Seifert fibred space with no singular
fibres. The surface 
$C$
is naturally contained in 
${\cal N}(K')$
and can be viewed as a circle bundle over a
disjoint union of simple closed curves in the
base surface of our Seifert fibration.

The boundary
$\partial {\cal N}(K')$
consists of tori. If a toral component 
of
$\partial {\cal N}(K')$
is compressible in 
$M$,
then it either bounds a solid torus in 
$M$ (which is disjoint from
${\cal N}(K')$,
because a solid torus can not contain an incompressible torus)
or it is contained in a 3-ball.
The latter possibility can not occur, because that would imply that 
some fibre from
${\cal N}(K')$
is homotopically trivial, contradicting the incompressibility
of patches of
$K'$.

In the former case we can extend the Seifert fibration over the 
solid torus. We can always do that because the fibres in the boundary
are homotopically non-trivial in
$M$.
Doing that for all compressible components of 
$\partial {\cal N}(K')$,
we obtain a new Seifert fibred space
$X$
which is embedded in
$M$,
contains
${\cal N}(K')$,
and must have a non-empty incompressible boundary
in
$M$.
If the 3-manifold
$X$
were closed, then it would have to be equal to the whole of
$M$
making it into a closed Seifert fibred space containing vertical
canonical tori. This is clearly a contradiction.
We should note at this point
that the construction up to now was based on the fact that the surface
$C$
is incompressible. The hypothesis that 
$C$
is in fact canonical will now come into play.

Assume now that some simple closed curve in the base surface
of
$X$,
that a toral component 
$C_0$
of
$C$
fibres over, is not boundary parallel in the base space. Then
we can find another simple closed curve in the base surface, such that 
$C_0$
can not be isotoped off the vertical torus above this curve by an isotopy
that is supported in
$X$.
Because
$\partial X$
is incompressible in 
$M$,
so is the vertical torus.
Hence, since
$C_0$
is canonical,  there exists a homotopy
$h:C_0\times I\rightarrow M$
that will move 
$C_0$
so that it is disjoint from
the vertical torus and from the boundary
$\partial X$.
The preimage of
$\partial X$
under
$h$
is a closed surface in the product
$C_0\times I$.
Using Dehn's lemma and the loop theorem we can modify the homotopy
so that 
$h^{-1}(\partial X)$
becomes a disjoint union of copies of
$C_0$
in the interior of the product 
$C_0\times I$.
Using the fact that two homotopic embeddings of an incompressible 
surface are isotopic (see corollary 5.5 in~\cite{wald}) implies
that we were isotoping 
$C_0$
over a boundary component of
$X$
and that
$C_0$
is parallel to it. But this contradicts our initial assumption.
So the collection of
simple closed curves in the base surface
of
$X$,
that 
$C$
fibres over,
must be boundary parallel.

The classification of surfaces implies that there exists
a family of disjoint annuli in the base surface of
the Seifert fibred space
$X$,
such that one boundary component of each annulus is a boundary
curve of the base surface. The union of the other boundary components
of our annuli are precisely the simple closed curves in 
the base surface of
$X$
that the
surface 
$C$
fibres over.

Even though the manifold
$X$
can contain
singular fibres, the above family of annuli contains no singular
points. This again follows from the fact that 
$C$
is canonical. If one of the annuli contained a singular point,
then we could easily construct an incompressible torus in
$X$ 
that could not be isotoped off 
$C$.

\noindent \textbf{Claim 3.} Each coefficient 
$k_i$
can either be 1 or 2.

Let's assume that 
$k_i$
is larger than 2 and let's look at some annular patch of 
$F_i$
with nonzero weight.
Since there are at least 3 copies of this patch 
contained in
$C$,
after using the projection of
$X$
onto its base surface, we can conclude that the images
of two adjacent copies of our patch belong to the same
boundary component of one of the annuli constructed above.
Furthermore there exists an arc 
$\alpha$,
properly embedded in this annulus, running between the two
adjacent copies of the projections of our patch. 
The arc
$\alpha$
chops off a disc in the annulus, bounded by
$\alpha$
on one side
and an arc
$\beta$
which is contained in the image of the projection
of
$C$.
Now we can use this disc to isotope
$\beta$
onto 
$\alpha$.
Since our annulus contains no singular points, we can
extend this isotopy, thus obtaining an isotopy
of
$C$,
which pushes the annulus over 
$\beta$
onto the annulus over 
$\alpha$.
By choosing the arc
$\alpha$
judiciously, we can make sure that this isotopy
reduces the weight of
$C$
since the patch of
$F_i$
we started with had positive weight.
But this is a contradiction because 
$C$
had minimal weight. So the claim follows.

The formula lemma 6.2 of~\cite{hass}, that gives a bound on the 
number 
of all
fundamental surfaces 
is bounded above by
$2^{70t^2}$.
Using a well known bound on the normal complexity 
of fundamental surfaces and the estimate 
$2\cdot20t^22^{7t}\leq 2^{10t^2}$
for 
$t\geq 2$,
we get the bound from the proposition.
\end{proof}

It should be noted that in the last proof we only needed the property that 
surfaces in 
$C$
were ``torus-canonical''. In
other words we only used the fact
that each torus in
$C$
can be isotoped off any incompressible 
torus in
$M$. 
This notion is very close to but not identical
with the notion of a canonical surface. In fact an essential 
torus is torus-canonical if and only if it is either canonical or it 
is parallel to a torus formed from a canonical annulus and an
annulus in
$\partial M$
(see proposition 4.2 in~\cite{neumann}).
This means that even if the 3-manifold
$M$
has non-empty pattern on its boundary, the proposition~\ref{prop:cantori}
still bounds the normal complexity of 
$P$-canonical tori in it.
Now we will describe 
how to construct canonical
annuli in triangulated 3-manifolds.

\begin{prop}
\label{prop:canannuli}
Let 
$M$
be an
irreducible 3-manifold with non-empty
boundary. 
Assume also that 
$M$
admits a possibly empty boundary pattern 
$P$
which is contained in the 1-skeleton of the triangulation 
$T$
of
$M$.
Let
$t$
be the number of tetrahedra in 
$T$.
Then every 
$P$-canonical annulus in 
$M$
can be isotoped into normal form so that it consists of
not more than 
$2^{80t^2}$
normal discs.
The same bound is valid if we replace a single
$P$-canonical annulus with a maximal collection of 
disjoint topologically non-parallel
$P$-canonical annuli. 
\end{prop}

\begin{proof}
Let 
$A$
be a maximal collection of disjoint topologically non-parallel
$P$-canonical annuli in
$M$.
Without loss of generality we can assume that the surface
$A$
is in normal form and that its weight is minimal. We can therefore
express it as a sum of fundamental surfaces
$A=k_1F_1+\cdots+k_nF_n$
where the coefficients 
$k_i$
are positive integers. Just like in the proof of proposition~\ref{prop:cantori}
we can show that the surfaces
$F_i$
are incompressible $P$-boundary incompressible and pure. By applying
lemma~\ref{lem:patches} in the same way as before, we can conclude
that 
$\chi(F_i)$
is zero. 

We now want to bound each integer 
$k_i$.
Like before we can express
$A$
as a sum
$A=k_iF_i+S$,
where 
$S$
is a pure normal incompressible $P$-boundary incompressible
surface in 
$M$
with zero Euler characteristic. Let's look at the polyhedron
$K=F_i\cup S$.
By using the procedure described in figure~\ref{fig:patch}
we can make sure that 
$K$
contains no adjacent trivial patches that are bounded by a 
simple closed curve from the intersection
$F_i\cap S$.
The same argument that justified claim 2 from the previous proof 
tells us that after this modification,
the polyhedron
$K$
contains no trivial patches that are disjoint from
$\partial M$.
It might however contain trivial patches that do meet the 
boundary of
$M$.

We start by looking at adjacent trivial patches in
$K$,
i.e. the ones that meet along an arc of intersection
$F_i\cap S$
(any such pair of trivial patches can not be contained in one 
of the 
surfaces 
$F_i$
or
$S$,
because both
$F_i$
and
$S$
are not a disc).
By doing regular alterations in the sum
$A=k_iF_i+S$
along each such arc from
$F_i\cap S$
and then performing an isotopy that is invariant on the pattern
$P$
we get a new normal surface 
$S'$
such that
$w(S)=w(S')$
(adjacent trivial patches must have the same weight because 
$A$
has minimal weight).
Furthermore the equality
$A=k_iF_i+S'$
holds for the isotoped surface
$A$
and the new polyhedron
$K'=F_i\cup S'$
contains no adjacent trivial patches. 
It follows directly that the surface
$A$
now contains no adjacent trivial patches either.
In fact it contains no trivial patches at all.
If it contained one, we could do the same construction
as in the second part of the proof of
lemma~\ref{lem:patches} 
in~\cite{mijatov1} to obtain a normal sum
$A=A'+B$
where the normal surface 
$B$
misses the 1-skeleton. This contradiction,
together with the fact that 
$F_i$
is connected, implies that 
all patches of 
$A$ 
are either annuli or discs of length four. The latter case means that
the patches are discs which are bounded by four arcs, two of them in
$(k_iF_i)\cap S'$
and the other two in
$\partial M$.

The former case, when the patches are annuli, is almost 
identical to the situation from the previous
proposition. Again the
regular neighbourhood 
$\mathcal{N}(K')$
supports a structure of a Seifert fibred space.
The surface 
$A$
is vertical in this Seifert fibration and is therefore
determined by a collection of properly embedded disjoint arcs
in the base surface of the fibration. All compressible
toral components of
$\partial \mathcal{N}(K')$,
that are disjoint from
$\partial M$,
have to bound solid tori in
$M$
which are disjoint from
$\mathcal{N}(K')$.
This follows by the same argument as in the proof
of~\ref{prop:cantori}. We can extend the fibration over all 
such solid tori. 
The annuli from
$\partial \mathcal{N}(K')-\partial M$
are incompressible simply because their core curves are generators
of the fundamental group of
$A$.
If such an annulus is 
$P$-boundary compressible, then it is parallel to a pure
annulus in
$\partial M$.
This parallelity region is topologically a solid torus.
We can extend the existing fibration of our manifold over it.
After we've adjoined all such parallelity regions to our submanifold,
we obtain a Seifert fibred space 
$X$
that contains 
$A$
as a vertical surface. 

Like in the proof of~\ref{prop:cantori} we can show that any 
vertical torus in 
$X$,
which is not 
parallel to a component of
$\partial M$,
can be isotoped off 
$A$ 
by an isotopy that is supported in
$X$.
Since the tori from
$\partial X$
are incompressible in 
$M$,
we have two possibilities. Either we are
able to isotope 
$A$
off all of the essential vertical tori in
$X$
or we are not.
In the former case the annuli of 
$A$
are boundary parallel in
$X$.
In other words their projection onto the base surface of
$X$
chops off discs from the base surface that contain 
no singular points.
We can therefore use the same technique 
that proved claim 3 from the proof of~\ref{prop:cantori}
to bound each 
coefficient
$k_i$.
In the latter case
there exists at least one annulus
$A_0$
in
$A$,
lying in a component
$X_0$
of
$X$
which is homeomorphic to 
one of the atoroidal Seifert fibred spaces (see figure 3 in~\cite{mijatov1}).
This is because
$X_0$
can not have any
essential vertical tori. Also, since the surface
$\partial X$
is incompressible in
$M$,
the boundary components of
$X_0$
that contain
$\partial A_0$
must be parallel to components of
$\partial M$.
The components of
$\partial X_0$
which are disjoint from
$A_0$
are not necessarily parallel into 
$\partial M$
and
can be canonical tori in
$M$.
In this topological setting it is not hard 
to prove, using our trivial patch
reduction techniques and theorem~\ref{thm:sum}, that 
$A_0$
is actually fundamental in
$M$
and that it therefore satisfies proposition~\ref{prop:canannuli}.

We can now assume that all patches of 
$A$
are discs of length four. This implies that
the regular neighbourhood
$\mathcal{N}(K')$
supports an 
$I$-bundle structure. 
The boundary of the manifold
$\mathcal{N}(K')$
is naturally divided into two bits. The horizontal part, which is contained
in
$\partial M$,
is just a
$\partial I$-bundle. The vertical part is a complement (in
$\partial \mathcal{N}(K')$)
of the horizontal boundary and it consists of properly embedded pure annuli in
$M$.
We now have the following claim.\\ 
\noindent \textbf{Claim.} If an annulus 
$V$
from the vertical boundary of
$\mathcal{N}(K')$
is compressible, then it bounds a pure submanifold of the form
$D^2\times I$
in
$M$,
such that 
$(D^2\times I)\cap \partial M=D^2\times \partial I$
and
$(D^2\times I)\cap V=\partial D^2\times I$.
This submanifold is disjoint from
$A$.

By compressing the annulus 
$V$
we obtain two pure properly embedded
discs 
$D_1'$
and
$D_2'$
in
$M$.
From the the definition of the pattern
$P$
and the fact that 
$M$
is irreducible it follows that 
$D_1'$
and
$D_2'$
are parallel to discs
$D_1$
and
$D_2$
in
$\partial M$.
The parallelity region between the discs
$D_i$
and
$D_i'$
is a pure 3-ball 
$B_i$
for
$i=1,2$.
Then there are two possibilities.

First, the balls 
$B_1$
and
$B_2$
are disjoint. Then the discs
$D_1$
and
$D_2$
are also disjoint and the 2-sphere 
$D_1\cup D_2\cup V$
bounds a 3-ball with the required 
properties. Also
$A$
has to be disjoint from this 3-ball because it 
consists of incompressible annuli. So the claim follows.

Second, the balls 
$B_1$
and
$B_2$
are nested. We can assume that
$B_1$
is contained in
$B_2$
and therefore
$D_1$
lies in
$\mathrm{int}(D_2)$.
The surface 
$A$
is disjoint from 
$V$
and hence does not meet the boundary circles
$\partial V=\partial D_1\cup \partial D_2$.
Since 
$A$
consists of incompressible annuli, both intersections
$A\cap D_1$
and
$A\cap D_2$
must be empty.
On the other hand it follows from the definition of
$V$
that there exists an embedded arc in
$\partial M$,
running from 
$\partial D_1$
to 
$\partial A$,
which is disjoint from 
$\partial D_2$.
In other words this means that 
$A\cap \mathrm{int}(D_2)$
can not be empty. This contradiction proves the 
claim.

We can now extend the 
$I$-bundle structure of
$\mathcal{N}(K')$
over every compressible vertical annulus in its 
boundary. If an incompressible vertical annulus is 
$P$-boundary compressible, then it has to be parallel
to a pure annulus in
$\partial M$.
The parallelity region, which is a solid torus, can be
used to construct a pure 
$\partial$-compression discs for some
annulus in 
$A$.
So this can not arise.
After adjoining all
the solid cylinders to
$\mathcal{N}(K')$,
we obtain a pure
$I$-bundle
$X$
in
$M$
whose vertical boundary consists of 
incompressible $P$-boundary incompressible
annuli. The horizontal boundary of
$X$
is a pure subsurface of
$\partial M$.
The collection of vertical annuli in 
$\partial X$
is not empty, since
otherwise 
$M$
would have to be homeomorphic to an
$I$-bundle over a closed (possibly non-orientable) surface with 
an empty boundary pattern. But such manifolds contain no canonical 
annuli.

We now apply the fact that the annuli in
$A$
are 
$P$-canonical.
It follows, in the same way as before, that 
$A$
can be isotoped off any annulus in
$X$,
which fibres over a simple closed curve in the base surface of
$X$,
by an isotopy that is
supported in 
$X$. 
This implies that the 
components of
$A$
have to be parallel to the vertical annuli from
$\partial X$.
So we can use the same strategy as in claim 3 
from the proof of~\ref{prop:cantori},
to show that the coefficient
$k_i$
can not be greater than
2. The only difference is that 
the weight reducing isotopy is supported
in a pure solid cylinder rather than in a solid torus.
This concludes the proof.
\end{proof}

It is probably true that the bound on normal complexity
of canonical surfaces given by propositions~\ref{prop:cantori}
and~\ref{prop:canannuli} is not the best possible. On the other
hand it is not hard to construct simple closed curves in a bounded
surface, none of which is boundary parallel, but whose ``normal'' sum is 
boundary parallel. The 
$S^1$-bundle over such a surface can be a component 
of the characteristic submanifold in some ambient 3-manifold.
The boundary torus of the 
$S^1$-bundle is certainly not going to be fundamental, but in this
setting it is clearly canonical.

It is worth pointing out that propositions~\ref{prop:cantori}
and~\ref{prop:canannuli} together imply the same bound on the normal 
complexity of the whole JSJ-system of surfaces in our 3-manifold
$M$.

The proof of proposition~\ref{prop:canannuli} clearly bounds the
normal complexity of every trivial $P$-essential annulus in
$M$. 
This is because such an annulus 
can also be isotoped off of any other $P$-essential 
annulus in
$M$.
But it is not hard to show that,
in an atoroidal 3-manifold where all $P$-essential annuli
are trivial and which is not a solid
torus, every trivial $P$-essential annulus
is isotopic to a fundamental annulus. 
First we isotope it so that it has minimal weight. 
If it is a sum of two surfaces, we can assume that they are
connected.
By applying lemma~\ref{lem:patches} we can conclude that
both of them have zero Euler characteristic, since they 
are both pure. By theorem~\ref{thm:sum}
they are both incompressible and 
$P$-boundary incompressible. 
At least one of the summands is bounded. It can not be a Moebius band
because our manifold is not a solid torus. So it is
a trivial $P$-essential annulus which can be isotoped off of our
original annulus by an isotopy that is invariant on the pattern.
If the other summand is closed, then these two annuli
have to be parallel.
This contradicts
our minimal weight assumption on the original trivial $P$-essential
annulus. If the other summand is also a trivial $P$-essential
annulus then, after removing all trivial patches in the normal sum,
we can conclude that both of the summands have to be parallel to
our original annulus. This contradiction proves our claim.

%% file: statement.tex
In order to
state theorem~\ref{thm:non-fibred}
we need to clarify what we mean by a fibre-free 3-manifold.
Before we define it we need to recall some standard terminology.
A surface bundle with an orientable fibre 
$S$
is just a mapping torus, i.e. a quotient
$S\times I/(x,0)\sim(\varphi(x),1)$,
for some orientation preserving surface automorphism
$\varphi\dv S\rightarrow S$. 
Since
$S$
is orientable
this construction gives an orientable 3-manifold. But for a non-orientable
surface
$R$
with a non-trivial two-sheeted covering
$S\rightarrow R$,
the mapping cylinder of the covering projection is an
orientable twisted 
$I$-bundle over 
$R$.
Gluing two such 
$I$-bundles
together along their horizontal boundaries
by an automorphism of
$S$
gives a 3-manifold
$N$ 
which is foliated by parallel copies of 
$S$
and the two copies of
$R$.
The leaves of this foliation are the
``fibres'' of a natural projection map
$N\rightarrow I$,
where the two copies of 
$R$
are the preimages of the endpoints of the interval
$I$.
Such a 3-manifold 
$N$
will be called a \textit{semi-bundle} (with fibre
$S$) over an interval 
$I$.
The surfaces 
$S$
and
$R$
can be either closed or bounded.

3-manifolds which are semi-bundles do sometimes arise in nature. 
In one dimension lower
for example, a 
Klein bottle is a semi-bundle with fibre
$S^1$,
since it splits as a union of two Moebius bands. 
The simplest
example in dimension 3 is a connected sum of two projective spaces
$\RR P^3\#\RR P^3$
where the fibre is a 2-sphere. 
On the other hand a semi-bundle structure can never arise in a knot
complement. This is because the boundary circles of the 
two non-orientable leaves would be disjoint curves in the
boundary torus and 
could therefore be capped off by the 
annuli they bound. This would then 
give a closed non-orientable
surface in
$S^3$.

Any semi-bundle over an interval can be viewed as a quotient
of the product
$S\times I$
with identifications
$(x,0)\sim(\alpha(x),0)$
and
$(x,1)\sim(\beta(x),1)$.
The homeomorphisms 
$\alpha$
and
$\beta$
are orientation reversing fixed point free involutions of the fibre
$S$.
Using this representation it is easy to see that
every semi-bundle admits a two-sheeted
covering space which is a surface bundle. The holonomy of the
surface bundle is the composition of
$\alpha$
and
$\beta$.
So for example every semi-bundle,
except the one mentioned above, is irreducible because the surface
bundle covering it is 
covered by $\RR^3$
which is irreducible. 
Recall that a \textit{Haken 3-manifold} is
an irreducible 3-manifolds with possibly empty incompressible
boundary which contains an injective surface different
from a disc or a 2-sphere.
Let's now define the class of 3-manifolds 
we will consider in theorem~\ref{thm:non-fibred}. 

\begin{defin}
A 3-manifold 
$M$
with a non-trivial JSJ-decomposition or
with non-empty boundary
is \textit{fibre-free} if 
none of the strongly simple pieces in
its JSJ-decomposition are homeomorphic to a surface bundle or
semi-bundle which contains
no closed injective surfaces. A closed atoroidal Haken 3-manifold
is \textit{fibre-free} if 
it is not homeomorphic to a closed
surface bundle over
$S^1$
or to a closed surface semi-bundle over
$I$.
\end{defin}

Notice that by our definition a fibre-free 3-manifold 
$M$
can contain
a bounded surface bundle or semi-bundle if it is a part of
the characteristic submanifold
$\Sigma$
or if it contains a closed injective surface. 
In other words if, say a surface bundle supports a Seifert fibration,
then the manifold 
$M$
is still fibre-free. 
Before we state theorem~\ref{thm:non-fibred}, we should remind ourselves that 
the exponent in the formula below, containing the exponential function
$e(x)=2^x$,
stands for the composition of the function with itself rather
than for multiplication.\\ 

\begin{thm}
\label{thm:non-fibred}
Let
$M$
be a fibre-free Haken 3-manifold. Let 
$P$
and
$Q$
be two triangulations of
$M$
that contain
$p$
and
$q$
tetrahedra respectively. Then there exists a sequence of
Pachner moves of length at most
$e^{2^{ap}}(p)+ e^{2^{aq}}(q)$
which transforms 
$P$
into a triangulation isomorphic to
$Q$.
The constant 
$a$
is bounded above by
$200$.
The homeomorphism of 
$M$,
that realizes this simplicial isomorphism, is supported in
the characteristic submanifold
$\Sigma$
of
$M$
and it does not permute the components of
$\partial M$.
\end{thm}

This theorem
gives a conceptually trivial algorithm for determining whether any 3-manifold
is homeomorphic to a complement of a given non-fibred knot in the 
3-sphere. Say we also wanted a simple 
procedure enabling us to determine whether 
any knot is the same as our given non-fibred knot. It is enough to 
establish whether their respective complements are homeomorphic and,
if they are, to determin if the homeomorphism maps the meridian
of one onto the meridian of the other.
If the boundary torus of this knot complement is
not contained in the characteristic submanifold, then 
the homeomorphism
from theorem~\ref{thm:non-fibred} equals 
the identity
on the boundary. If on the other hand the bounding torus is contained in
$\Sigma$,
then we first make sure that the simplicial structures on the boundary of
both knot complements coincide. 
It will be clear from the proof of theorem~\ref{thm:non-fibred}
that this makes the homeomorphism equal
to the identity on the boundary torus.
So in this way, using theorem~\ref{thm:non-fibred}, we can solve the 
recognition problem for any non-fibred knot.

The proof of theorem~\ref{thm:non-fibred} starts by subdividing
the original triangulation of 
$M$
so that the characteristic submanifold 
$\Sigma$
is triangulated by a subcomplex of the subdivision. We then look at the
strongly simple pieces of the JSJ-decomposition. We will simplify
their triangulations using the canonical hierarchy which
is described in section~\ref{sec:can}. 
The ``canonical'' triangulation is
obtained from the canonical
hierarchy by applying theorem 1.2 from~\cite{mijatov}.
So a triangulation
of every strongly simple piece 
will impose a simplicial
structure in the boundary of both
Seifert fibred and 
$I$-bundle components
of
$\Sigma$.
The former can then be simplified by theorem 3.1 of~\cite{mijatov1}, 
and the latter can be dealt with using techniques from
subsection 6.2 of~\cite{mijatov1} and
subsection~\ref{subsec:norcom}.

%% file: canonical.tex
Let $M$
be a 3-manifold satisfying the hypothesis of 
theorem~\ref{thm:non-fibred}.
In this section we are going to describe the
canonical two-dimensional object lying inside 
$M$,
that will give rise to an intermediate triangulation 
which will later be used to bridge the 
gap between any two triangulations of
$M$.
Let's start by introducing some 
standard terminology. 
Let 
$S$
be an incompressible surface 
contained in the 3-manifold 
$M$.
The operation of 
\textit{cutting $M$ along $S$}
results in a 3-manifold
$M_S$,
that is just a complement 
(in
$M$)
of
the interior of the regular 
neighbourhood 
$\cal{N}(S)$.

\begin{defin}
A \textit{partial hierarchy} for a Haken
3-manifold
$M$
is a sequence of 3-manifolds
$M_1, \ldots, M_n$,
where
$M_1$
equals
$M$,
and 
$M_{i+1}$
is obtained from 
$M_i$
by a cutting along an orientable, incompressible,
properly embedded surface in 
$M_i$,
no component of which is a 2-sphere.
A \textit{hierarchy} is a partial hierarchy with
the property that 
$M_n$
is a collection of 3-balls.
We shall denote (partial) hierarchies as follows
$$M_1 \stackrel{S_1}{\longrightarrow}
M_2 \stackrel{S_2}{\longrightarrow} \cdots
M_{n-1}\stackrel{S_{n-1}}{\longrightarrow} M_n,$$
where 
$S_i$
is the surface in 
$M_i$
that we cut along. 
\end{defin}

It is well known
that every
Haken 3-manifold possesses a hierarchy. Since we are going
to construct a hierarchy with some additional properties 
we will be, among other things, reproving this 
result.
Another classical fact about hierarchies is contained in
lemma~\ref{lem:inco}. Since 
it represents a key step 
in the construction of the canonical hierarchy,
we will include its proof.\\ 

\begin{lem}
\label{lem:inco}
Let
$M$ 
be a compact orientable irreducible 3-manifold with
(possibly empty) boundary. Let 
$M_1 \stackrel{S_1}{\longrightarrow}
M_2 \stackrel{S_2}{\longrightarrow} \cdots
M_{k-1}\stackrel{S_{k-1}}{\longrightarrow} M_k$
be a partial hierarchy for 
$M$
and let 
$N$
be a regular neighbourhood of the
union of all the surfaces in the 
hierarchy, i.e. 
$N=\mathcal{N}(\partial M\cup S_{1} \cup\ldots\cup S_{k-1})$.
Then the following holds:
\vspace{-1mm}
\begin{itemize}
\item[1.] The surface
          $\partial N - \partial M$
          is incompressible in 
          $N$.
\vspace{-2mm}
\item[2.] Compressing the surface 
          $\partial N - \partial M$
          into the 3-manifold
          $M - \mathrm{int}(N)$
          as much as possible gives a disjoint union
          of closed connected separating surfaces, each
          of which is either incompressible in
          $M$
          or is a 2-sphere bounding a 3-ball in
          $M$. 
\end{itemize}
\end{lem}

\begin{proof}
Let's start by observing that the (partial)
hierarchy for
$M$
from the lemma gives a partial hierarchy for
$N$,
where the manifold we end up with after cutting along
surfaces
$S_1,\ldots, S_{k-1}$
is homeomorphic to
$(\partial N - \partial M)\times I$.
Now we proceed by induction on
$k$. 
For
$k$ 
equals 2 we have a single incompressible surface 
$S_1$.
We can assume that a compression disc
$D$
for
$\partial N - \partial M$
(in $N$)
intersects 
$S_1$
in simple closed curves only and that it is 
disjoint from 
$\partial M$.
The innermost simple closed curve in
$D$
bounds a disc  
$D_1$
in
$S_1$
(since 
$S_1$
is incompressible) 
and a disc
$D'$
in 
$D$.
The 2-sphere
$D'\cup D_1$
bounds a 3-ball in
$M$,
since 
$M$
is irreducible. This 3-ball might not
be contained in 
$N$,
but its boundary certainly is.
So by isotoping over the 3-ball, we obtain a new disc 
(again denoted by 
$D$),
lying inside
$N$
such that
the number of components of the intersection
$D\cap S_1$
is reduced by at least one.
Repeating this procedure, we can make 
$D$
disjoint from
$S_1$.
This implies that 
$D$
in fact is a compression disc for the surface
$\partial N - \partial M$
in the product
$(\partial N - \partial M)\times I$.
So its boundary has to bound a disc in
$(\partial N - \partial M)$.
The inductive step is proved in the similar fashion
by first making the compression disc disjoint from the 
surface
$S_1$,
and then applying the induction hypothesis to the
partial hierarchy in
$M_{S_1}$. 
This proves 1. 

The closed surface
$\partial N - \partial M$
is separating in
$M$
and can only compress into
$M - \textrm{int}(N)$.
Furthermore, these new compression discs
can be 
viewed as the continuation of the existing
partial hierarchy. 
So by applying the first part of the lemma to
this extended hierarchy, we get that each component of the 
compressed surface is either incompressible both 
in the regular neighbourhood of the extended
hierarchy as well as in its complement or it is a 2-sphere.
In both cases part 2 of the lemma follows. 
\end{proof}

The canonical two-dimensional polyhedron
lying inside of the hyperbolic pieces of 
$M$
can now be topologically,
but not yet algorithmically, described as follows. 
It will consist of the union of surfaces in the
canonical hierarchy which is given by the following
three stage procedure: 
\newcounter{c}
\begin{list}{\arabic{c}.}{\usecounter{c}\parsep 0mm\topsep 1mm}
\item 
      Let 
      $S_1$
      be the first surface 
      in the hierarchy. It's components
      will be defined recursively.
      We start by adjoining
      all canonical tori and annuli
      from the JSJ-system of
      $M$.
      In other words at this point 
      $S_1$
      consists of the surfaces
      $\partial \Sigma-\partial M$,
      where 
      $\Sigma$ 
      is the characteristic submanifold of
      $M$.
      We will adjoin new closed connected surfaces to
      $S_1$
      in the components of the JSJ-decomposition of
      $M$
      which are 
      not contained in 
      $\Sigma$
      and are neither homeomorphic to a solid torus nor to a
      $(\mathrm{torus})\times I$.
      Boundary of such a piece
      might or might not be incompressible in
      $M$.
      For each incompressible boundary component of such
      a component, we add a parallel
      copy of it to our surface
      $S_1$
      (this makes sure that all the complementary pieces 
      we need to work on in step 2 are pure).
      In each component of the complement 
      $M-\mathrm{int}(\Sigma)$ 
      defined above,
      the rest of the
      surface 
      $S_1$
      will
      consist of 
      disjoint closed connected
      orientable incompressible surfaces. If two components
      of
      $S_1$
      are parallel, then no other connected subsurface in 
      $S_1$
      will be parallel to them.
      They are defined by the following recursion. 
      \begin{itemize}
      \item[(1a)] If the complementary piece of the surface we've defined
                 so far contains no 
                 canonical annuli, then we proceed
                 by looking for a closed, 
                 injective surface that is not boundary
                 parallel (in that piece), and that has the largest
                 Euler characteristic (i.e. smallest genus) out of all
                 such surfaces in our piece. If a surface 
                 like that does not exist, then 
                 $S_1$
                 gets no new components in the piece we are studying.
                 If it does exist, then we adjoin to
                 $S_1$
                 a boundary of its regular neighbourhood in our piece.

      \item[(1b)] If, on the other hand, our submanifold contains a 
                 canonical annulus,
                 we look at the union of all 
                 canonical annuli with the boundary of the piece.
                 Compressing the boundary of the regular neighbourhood of that 
                 two-dimensional polyhedron as 
                 far as we can, gives a disjoint union
                 of closed separating orientable surfaces. 
                 If this union contains a surface
                 that is neither
                 homeomorphic to 
                 $S^2$,
                 nor parallel to any of the boundary components of the 
                 piece, 
                 then we adjoin to 
                 $S_1$
                 a boundary of a regular neighbourhood of
                 a closed injective
                 surface (in our piece), which is not boundary
                 parallel, and that has 
                 the largest negative Euler characteristic 
                 among all such surfaces.
                 If there aren't any such surfaces, we do nothing.
      \end{itemize}
      The union of all the surfaces in
      $M$
      obtained in this way, defines the
      surface 
      $S_1$
      in our hierarchy.

\item Let's look at a piece of the cut-open manifold
      $M_{S_1}$
      that is not homeomorphic to an 
      $I$-bundle over a surface and is not contained in the
      characteristic submanifold of 
      $M$
      (notice that all such pieces are pure).
      Such a piece either contains a canonical annulus or it doesn't. 
      If it does, then take the surface
      $S_2$
      in that piece
      to be a union of two parallel copies of each 
      canonical annulus in it.
      If it doesn't, then let 
      $S_2$
      be the horizontal boundary of a regular
      neighbourhood of a connected bounded 
      two-sided incompressible surface with the
      largest Euler characteristic. 

\item The pieces of 
      $M$ 
      after the first two steps are
      either homeomorphic to
      compression bodies, whose ``negative'' boundary
      $\partial_-$
      can be empty
      (i.e. handlebodies) or disconnected, or to  
      $I$-bundles over (perhaps non-orientable) closed surfaces.
      All the subsequent surfaces in the canonical hierarchy
      are going to be annuli and discs and are going to be defined
      recursively using boundary patterns on the pieces.
\end{list}

The above construction both raises a number of questions 
and also requires some explanation. 
The recursion defining  
the surfaces 
$S_1$
must stop by the Kneser-Haken
finiteness theorem (see~\cite{jaco1} and~\cite{hempel}). Moreover 
it follows that
$M$
can contain at most 
$8t+\beta_1(M;\ZZ)+\beta_1(M;\ZZ_2)$
closed connected incompressible surfaces that are non-parallel
($t$
is
the number of tetrahedra in $T$).
Therefore
the surface
$S_1$
can contain at most
$2\cdot 20t$
components since
the Betti numbers
$\beta_1(M;\ZZ)$
and
$\beta_1(M;\ZZ_2)$
are smaller than 
$6t$. 

It is clear from the definition of 
$S_1$
that all components of the cut-open manifold obtained
at each stage of the recursion defining 
$S_1$,
which are disjoint from the characteristic submanifold
$\Sigma$,
are atoroidal. So the set of all surfaces in the
JSJ-system of each such 3-dimensional piece is either 
empty or it contains only annuli. 
All components of the 3-manifold
$M_{S_1}$
have non-empty boundary. It is a very well known fact
(see~\cite{hempel}) that a bounded compact
irreducible 3-manifold, which is not a 3-ball, contains
an orientable (and hence two-sided) 
non-separating properly embedded incompressible
surface with boundary. So even in the components of
$M_{S_1}$
that contain no canonical annuli, we can still find and fix one such surface
with the largest Euler characteristic.
Therefore we can carry out step 2
as well.

In order to determine what the pieces of the complement 
look like topologically after stage 2 we need the following 
definition. A \textit{compression body}
is a connected 
orientable 3-manifold obtained from a (possibly disconnected)
closed orientable surface
$S$
by attaching 1-handles to
the surface
$S\times \{1\}$
in the boundary of the product
$S\times I$.
The ``negative'' boundary 
$\partial_-$
of the compression body is precisely the incompressible part
of its boundary, which equals 
$S\times \{0\}$.
The rest of the boundary is called ``positive'' and is denoted
by
$\partial_+$.
It is convenient to regard handlebodies as compression bodies
with empty ``negative'' boundary
$\partial_-$.

We now have to examine 
the boundary of the regular neighbourhood
$N=\mathcal{N}(\partial M\cup S_1\cup S_2).$
The boundary of 
$N$
is a disjoint union of closed surfaces some of which are parallel to
the components of the surface
$\partial M \cup S_1$
and some of which are not. The former ones are incompressible in
$M$.
The latter ones are incompressible only in
$N$.
Lemma~\ref{lem:inco} also tells us that these
surfaces either compress all the way, and therefore
bound handlebodies, or they become closed incompressible   
surfaces in 
$M$ 
after a few compressions.
The second alternative implies that they have to be parallel 
(in
$M$)
to components of 
$\partial M\cup S_1$.
This follows from the construction of
$S_1$,
since all
closed orientable incompressible surfaces 
in
any component of the 3-manifold
$M_{S_1}$
have to be boundary parallel.
Putting all this together we can conclude that each
piece of the cut-open manifold
$M-\mathrm{int}(N)$ 
is homeomorphic to one of the following:
a handlebody, 
a compression body, with a possibly disconnected
``negative'' boundary
$\partial_-$, or an
$I$-bundle over a (possibly non-orientable) closed 
incompressible surface.

All surfaces in the third stage of our hierarchy will be annuli
and discs. The topological information coming from the manifold
$M$
will be reflected through the boundary patterns
(as described in
subsection~\ref{subsec:Pcan}) on the boundaries of the
pieces. 

There is a very natural way of obtaining boundary
pattern from any partial hierarchy of
$M$.
Let  
$M_1 \stackrel{S_1}{\longrightarrow}
M_2 \stackrel{S_2}{\longrightarrow} \cdots
\stackrel{S_{k-1}}{\longrightarrow} M_k$
be such a hierarchy and let
$K$
denote a two-dimensional polyhedron which is a union of 
$\partial M$
with all the surfaces from this hierarchy.
The singular locus 
$S(K)$
of the polyhedron
$K$
consists of all the points in 
$K$
that do not have neighbourhoods homeomorphic to
discs (i.e. the points whose links
in
$K$
are not circles).
Let
$Q$
be a closure (in
$M$)
of
a component 
from
$M-K$.
In our setting 
$Q$
is always an embedded submanifold of
$M$.
It also inherits the boundary pattern from the hierarchy
in the following way:
$P=\partial Q\cap S(K)$.
We should note that 
$P$
consists of three-valent graphs and disjoint simple closed curves,
and that the surface
$\partial Q- P$
is by definition incompressible in
$Q$.

Let now 
$Q$
be a complementary piece 
the canonical hierarchy after step 2.
Observe that 
the pattern on 
$\partial Q$
consists of a (possibly empty) collection
of disjoint simple closed curves.
All the subsequent surfaces in the canonical hierarchy
which are contained in
$Q$,
are going to be annuli and discs. They will
be defined 
using the information coming from the pattern.
We will adjoin them in such a way so that the closures 
of the complementary
components at any stage of the hierarchy, are embedded 
submanifolds. 

Before we proceed to describe the next step
in the canonical hierarchy, we need to impose some pattern
on the pure pieces of the form
$(\mathrm{torus})\times I$
which are not contained in the characteristic submanifold
$\Sigma$. 
Such submanifolds arise from canonical tori in
$M$
that are separating Seifert fibred pieces with
non-matching fibrations.
So on each
boundary component of
$(\mathrm{torus})\times I$,
we take our pattern to be the regular fibre of the Seifert fibration
on that side. Since the fibrations do no match (if they did, the
defining torus of the piece would not be canonical) this pattern 
admits no non-trivial $P$-essential annuli in our piece.
The third step in the construction of the canonical hierarchy
will be carried out in three substages as follows:

\newcounter{d}
\begin{list}{(3\alph{d})}{\usecounter{d}\parsep 0mm\topsep 1mm}
\item       Fix an ordering of the 3-manifold pieces
            in the complement 
            $M-(S_1\cup S_2)$
            that are not contained in the characteristic
            submanifold
            $\Sigma$.
            Go through all the pieces one by one, 
            respecting this ordering,
            and each time add a disjoint maximal
            collection of $P$-canonical annuli 
            that exist in that piece. 
            Every time we adjoin new
            $P$-canonical annuli, the pattern on the boundaries 
            of the neighbouring pieces acquires some
            additional simple closed curves.
            Also if a $P$-canonical annulus does not separate
            the piece it lies in, we add two copies of it.
            When we are done with all the components of 
            $M-(S_1\cup S_2)$,
            we fix a new ordering on the pieces of
            $M-(S_1\cup S_2\cup S)$,
            where 
            $S$
            equals the union of all $P$-canonical annuli
            we've added in so far. We repeat the above
            procedure on this new ordering. 
            We reiterate the whole process until we reach
            an ordering of the complementary components
            with boundary patterns admitting no 
            $P$-canonical annuli in any of the pieces.

\item       At this stage of the hierarchy all
            the complementary pieces which contain
            non-trivial 
            $P$-essential annuli are 
            $I$-bundles with pure
            horizontal boundaries (see theorem~\ref{thm:bundle}).
            In this step we are going to simplify
            proper compression bodies and 
            $I$-bundles over closed surfaces down to handlebodies
            by inserting the spanning annuli.
            First we fix a complementary piece that is a compression body.
            We adjoin each trivial $P$-essential annulus in our
            piece with the following property:
            the annulus in the boundary
            of the piece that it is parallel to 
            contains precisely two pure sub-annuli 
            and no simple closed curve components 
            of the pattern (this technical requirement
            is here to insure that normal surface theory
            can be applied later on and is not otherwise essential).
            We can do this
            step by step 
            by adjoining one such annulus at a time. This does not
            alter the topology of our piece.
            Now we
            adjoin an incompressible spanning annulus that has minimal 
            intersection with the pattern.
            Spanning in the
            context of compression bodies means that the bounding circles 
            of the annulus 
            lie in distinct boundary 
            components of the compression body. 
            If the chosen spanning annulus is not separating,
            we take two parallel copies of it.
            We repeat this procedure for each component of the
            ``negative'' boundary
            $\partial_-$.
            We now run
            step (3a) again. 
            If there are compression bodies left we 
            repeat (3b) on one of them, otherwise we do (3b) on an
            $I$-bundle over closed surface. Once all complementary pieces
            are handlebodies we move to the next step.

\item       Now all complementary pieces are 
            handlebodies and 
            none of them contain a $P$-canonical annulus. 
            Fix a complementary
            piece that contains no 
            non-trivial $P$-essential annuli and which is neither a 3-ball
            nor a solid torus with a pure injective annulus in its 
            boundary.
            If no such piece exists we proceed to the next step.
            Otherwise we adjoin all trivial $P$-essential annuli in our
            piece that have the same property as the ones described
            in step (3b). We adjoin them using the same method as
            in (3b).
            Then
            we choose a compression 
            disc in our piece which 
            has minimal
            intersection number with the boundary 
            pattern among all compression discs in our handlebody. 
            If it is not separating, we take two copies of it. 
            Notice that the compressed handlebody, which 
            is not necessarily connected, satisfies the same conditions
            as the original piece we've just compressed. So 
            we carry on with the compressions until the original
            piece becomes a union of 3-balls. 
            Now we run step (3a) again.
            The situation now is precisely as it was at the beginning of
            (3c), so we can repeat it.

\end{list} 

The description of step 3 requires some explanation. Step (3a)
is there to eliminate the non-trivial $P$-essential annuli in the
complementary pieces. This is crucial for the algorithmic construction
of the surfaces in the hierarchy, because such annuli make
normal surface theory impossible to apply. In step (3b)
we make sure that all complementary pieces are handlebodies. 
Step (3c) is there to compress these handlebodies, which
are different from a solid torus with a pure annulus
in its boundary, down to
3-balls.
 
Clearly if no complementary
piece after (3a) is homeomorphic 
to a proper compression body or to an
$I$-bundle over a closed surface,
we do nothing in (3b). If on the other hand one such piece exists,
then it contains no non-trivial $P$-essential
annulus. This will follow from theorem~\ref{thm:bundle} which is 
applicable since
all the complementary pieces we are looking at in step 3 are atoroidal. 
Therefore a compression body,
containing a non-trivial $P$-essential annulus but not a $P$-canonical
one, would have to be homeomorphic to either an atoroidal Seifert fibred 
space which is not a solid torus
or to an I-bundle. The first possibility can not occur since 
all such Seifert fibred spaces (see figure 3 in~\cite{mijatov1})
have incompressible boundary.
Also a proper compression body is not homeomorphic to an
$I$-bundle over a compact surface because it has a compressible
boundary component and at least one incompressible 
one. Since both possibilities lead to contradiction,
no compression body can contain a non-trivial $P$-essential annulus
after step (3a). Notice also that the complement of the spanning 
annulus in a compression body is again a compression body 
whose ``negative'' boundary has fewer components than that 
of the original piece.

If there exists an
$I$-bundle over a closed surface which 
contains a non-trivial $P$-essential 
annulus, 
then by theorem~\ref{thm:bundle} it
would have to be homeomorphic to 
an
$I$-bundle over a
compact surface with pattern lying in the annuli that
fibre over the boundary circles of the base surface.
So our piece is actually
pure because any
$I$-bundle over a bounded surface has compressible boundary. This
means that the boundary of our 
$I$-bundle is contained in the surface
$S_1$.
Since not more than two components of
$S_1$
are parallel, our piece can either be adjacent to another 
$I$-bundle or
the complementary piece lying on the other side of the boundary 
is a handlebody or a compression body. 
The first scenario would make
$M$
into a closed atoroidal surface bundle (respectively semi-bundle)
over 
$S^1$
(respectively
$I$).
This contradicts our initial hypothesis on
$M$.
The second scenario with the compression body can not take place,
because step (3b) does not touch 
$I$-bundles while compression bodies 
are still around. But the adjacent piece
can also not be a handlebody since all components of
$S_1$
are incompressible in
$M$.
All this implies that, when there are no more
compression bodies in the complement of the hierarchy, each 
$I$-bundle over a closed surface contains no non-trivial
$P$-essential annuli.
After we adjoin a spanning annulus our 
$I$-bundle becomes a handlebody.

In steps (3b) and (3c) we also have to adjoin
certain trivial $P$-essential annuli to the complementary pieces. 
We add in the ones that are parallel to the annuli in the 
boundary of the piece which contain no simple closed 
curves of the pattern and do contain precisely two pure 
sub-annuli. This rule makes sure that no two of the 
trivial $P$-essential annuli we add are topologically 
parallel (i.e. parallel disregarding the pattern).
This process can therefore not change the
topology of the piece. 
It only creates more of the solid tori,
in the complement of the hierarchy, that contain injective pure
annuli in their boundaries. 
The reason why we have
to adjoin the trivial $P$-essential annuli before compressing,
is because we need to make sure that any surface 
$F$
we are trying to
find, be it a disc or an annulus, 
is isotopic, by an isotopy that is invariant on the pattern,
to a surface obtained from 
$F$
by twisting along a trivial $P$-essential annulus. 
It is clear that after we adjoin the trivial annuli this can
be achieved.

Now that we've fully described the canonical hierarchy, we need to
show that step 3 (i.e.
(3a) and hence most other substeps) 
does not run forever and that, when it terminates, 
the complementary pieces of the hierarchy we get are
solid tori with pure annuli in their boundaries and 3-balls.
To do that we have to understand
how incompressible (but not
$\partial$-incompressible) annuli lie in handlebodies and compression bodies.
In fact we can concentrate only on the annuli that are not 
$\partial$-parallel in the piece we adjoin them to, because
the number of the boundary parallel ones 
is easily controlled by the complexity of the pattern.

It follows directly from the definition that a $P$-canonical annulus 
is not boundary parallel in the complementary piece we adjoin it to.
Let's start by looking at a handlebody component 
$H$
of
$M-\mathrm{int}(\mathcal{N}(S_1\cup S_2))$.
The process described above creates a sequence of pure incompressible
properly embedded annuli in 
$H$. 
Moreover no two annuli can be parallel because that would make
one of them 
$\partial$-parallel in a 
submanifold of
$H$.
So what we have is a disjoint collection of properly embedded 
non-parallel annuli in the handlebody 
$H$
(here we are only taking into account one of the two
parallel copies of $P$-canonical 
annuli that we sometimes had to adjoin, because the annulus
was not separating)
none of which is 
$\partial$-parallel.

If 
$g$
is the genus of handlebody 
$H$,
then there are at most 
$6g$
such annuli. 
If
$\mathcal{A}$
is the union of at least
$6g$
such annuli in
$H$,
we can perform a sequence of
$\partial$-compressions to 
$\mathcal{A}$
until we end up with a 
$\partial$-incompressible surface.
Since the only such surfaces in handlebodies are
discs, we get a collection of at least 
$6g$
disjoint discs in
$H$.
Notice that each annulus in
$\mathcal{A}$
is compressed only once and that the disc it yields, after the
$\partial$-compression, is a compression disc for 
$\partial H$.
This is because none of the components of
$\mathcal{A}$
are 
$\partial$-parallel.
But there are at most
$3g-3$
disjoint non-parallel compression discs in a handlebody of 
genus
$g$.
So at least three of our discs are parallel. 
It is easy to see that, if we reconstruct
the annuli corresponding to these parallel discs by reversing
the
$\partial$-compressions, at least two of the annuli 
are going to be parallel in
$H$.

The same proof tells us that in a compression body 
$H$
we can not have
more than
$6g(\partial_+H)$
incompressible
$\partial$-compressible annuli which are nether
$\partial$-parallel nor parallel to each other
(a compression body contains at most 
$3g(\partial_+H)-3$
disjoint non-parallel compression disc).
On the other hand compression bodies do contain incompressible
$\partial$-incompressible annuli. 
We can always choose the handle
structure of the compression body
$H$
so that our family of annuli is vertical
in the product structure of the complement of the 1-handles. 
This can be seen as follows. Choose a family of compression
discs in
$H$
which cuts it down to
$(\partial_-H)\times I$
and which intersects our annuli in the minimal number of arcs.
Since all the annuli 
are 
$\partial$-incompressible,
all the arcs of intersection must be inessential. Consider
an outermost arc separating off a disc 
$D$
in one of the annuli. The arc lies in one of the
compression disc and divides it into subdiscs
$D'$
and
$D''$.
Remove
$D'$,
say, from the compression disc and replace it with
$D$.
By choosing 
$D'$
appropriately we may ensure that the  family 
of the compression discs we obtain after
this operation still cuts 
$H$
down to 
$(\partial_-H)\times I$.
Clearly this family of discs intersects our annuli in 
fewer arcs which is a contradiction.

This implies that there are no more than 
$(3g(\partial_-H)-3)+2(g(\partial_+H)-g(\partial_-H))$
non-parallel incompressible 
$\partial$-incompressible annuli in 
$H$.
The first summand bounds the number of disjoint non-parallel
simple closed curves in the surface
$\partial_- H$.
The second summand accounts for the fact that two vertical annuli
which are parallel in
$\partial_- H\times I$
are not necessarily parallel in
$H$. 
Since each 1-handle of 
$H$
produces two discs in the boundary of
$(\partial_- H)\times I$
we get the above bound. Putting all these numbers together
we can conclude that each compression body can contain
not more than
$9g(\partial H)$
disjoint non-parallel incompressible annuli which are not 
$\partial$-parallel.

We are going to give a bound on the genus of any closed surface,
bounding a compression body or a handlebody from
$M-\mathrm{int}(\mathcal{N}(S_1\cup S_2))$
(see corollary~\ref{cor:topcomp}).
This implies that there are only finitely many compression
discs in step (3c). The number of annuli in steps (3a), (3b) and
(3c), which are not boundary parallel, is, by  
the
discussion in the previous
paragraphs, controlled by the genus of these surfaces.
Each trivial $P$-essential annulus that features in steps (3b) and (3c)
is always parallel to an annulus in the boundary of the piece which
contains some pattern. This pattern is such that it must contain
at least one simple closed curve, which is embedded in it, that
is a boundary component of some surface in the hierarchy that
was adjoined 
to the neighbouring 
piece at some earlier instance
(notice that this simple closed curve must contain at 
least one vertex of the pattern). 
Corollary 5.8 in~\cite{wang} gives us control over the negative
Euler characteristic of non-annular components of 
$S_2$.
The number~\label{bound} of boundary components of annular pieces of 
$S_2$
is bounded by 
$20\cdot18 t$
where 
$t$
is the number of tetrahedra in 
$M$.
This follows from the Kneser-Haken finiteness theorem and
the fact that all canonical annuli lie in the complement
of the amalgams (see the proof of corollary~\ref{cor:topcomp}).
So we have complete control over the number
of boundary components of the surface
$S_2$.
Also when we adjoin trivial $P$-essential annuli, the 
pattern that they generate can not contribute to 
the birth of new trivial $P$-essential annuli. Therefore,
since we can bound the number of $P$-canonical annuli,
we can also bound the number of the trivial
$P$-essential ones that come out of (3b) and (3c).
All this implies  that the procedure described by step 3 in the 
construction of the canonical hierarchy must terminate.

We now need to make sure that all the embedded submanifolds in
the complement of the whole canonical hierarchy, after 
we've carried out step 3, are indeed 3-balls and 
tori that contain pure annuli in their boundaries. Since step
(3a) can not run forever under any
circumstances, we can assume 
that steps (3b) and (3c) have been implemented as well.
This means that all the complementary pieces are handlebodies
that contain no $P$-canonical annuli and 
3-balls. 

Assume now that there exists a
piece 
$Q$
which contains a non-trivial $P$-essential annulus. Then,
by theorem~\ref{thm:bundle}, it has to be an
$I$-bundle over a bounded surface 
of negative Euler characteristic. 
The pattern is contained in the vertical boundary
and the horizontal boundary 
is pure.
Now we look at the complementary piece which is on the other side of 
the horizontal boundary
of
$Q$. 
The horizontal
boundary is incompressible in all pieces adjacent to 
$Q$. 
These pieces are therefore not 3-balls. Furthermore
such a piece must contain a non-trivial $P$-essential annulus, because if it
didn't, its boundary would have been compressed down to a 2-sphere
by step (3c). Since the adjacent piece doesn't contain any $P$-canonical annuli,
we can apply theorem~\ref{thm:bundle} to conclude that 
it is an
$I$-bundle over a bounded surface as well. Moreover the horizontal
boundary of that 
$I$-bundle is pure. It 
is either connected or it has two components, depending on whether
the base surface is non-orientable or not. In the first case the horizontal
boundary of the 
$I$-bundle adjacent to
$Q$
has to coincide with a component of the horizontal boundary
of 
$Q$. 
In the second case a component of the
horizontal boundary 
matches a component of the horizontal 
boundary of
$Q$.
If 
$Q$
and its neighbour meet along all of their horizontal boundaries,
then their union is an embedded submanifold
$N$
in
$M$
that is a surface bundle
or semi-bundle 
and the fibre is a bounded surface.
If not
we can extend the 
$I$-bundle structure over the union of 
$Q$
and its neighbouring piece. We can repeat the same argument
for the enlarged 
$I$-bundle and the pieces adjacent to it.
After finitely many repetitions we must arrive at the surface
(semi-)bundle situation described above, since
we are in the complement of
$\Sigma$.
The following lemma will
give the final contradiction with our initial hypothesis on the 
3-manifold
$M$.

\begin{lem}
\label{lem:contra}
The submanifold
$N$
constructed above 
is a bounded strongly simple piece of the 
JSJ-decomposition of the 3-manifold
$M$.
In particular 
$N$
is both atoroidal and an-annular.
Recall that the surfaces
$S_1$
and
$S_2$
constitute the first and the second step in the construction
of the canonical hierarchy.
We also have that
$\partial N-\partial M$
equals
$N\cap S_1$
and that it
consists of canonical tori from
$M$.
The intersection 
$N\cap S_2$
is either a union of two parallel components of the surface 
$S_2$
(if 
$N$
is s surface bundle over 
$S^1$)
or a single component of
$S_2$
(if 
$N$
is s surface semi-bundle over 
$I$).
\end{lem}

\begin{proof}
Since 
$M$
is fibre-free and the submanifold
$N$
lies in the complement of
$\Sigma$,
the surface fibre of 
$N$ 
can not be closed. In other words the boundary of
$N$
can not be empty.
First we need to show that the components
of
$\partial N$
are incompressible tori in
$M$.
Let 
$S$
be one such torus. Assume that it is compressible. Then it 
either bounds a solid torus in
$M$
or is contained in a 3-ball in
$M$.
But 
$S$
contains a boundary circle of the fibre of 
$N$. 
Since this surface injects in
$M$,
the torus 
$S$
can not be contained in a 3-ball nor can it bound a solid
torus with the slope of the meridian disc equalling the slope 
of the fibre of
$N$.

Let 
$F$
be a connected surface in the canonical hierarchy that contains
one of the fibres of the surface bundle
$N$. 
It is clear that 
$F$
has to be either a component of the closed surface 
$S_1$
or a non-annular component of the bounded surface 
$S_2$.
In both cases we are going to 
consider how the surface 
$F$
interacts with the solid torus bounded by
$S$
(we can assume without loss of generality that the solid 
torus is disjoint from
$\mathrm{int}(N)$).
We have two possibilities. The surface 
$F$
is either disjoint from its interior, or
it intersects it. In the former case it has to
contain another surface fibre of
$N$.
This leads directly into contradiction because it
renders
$F$
compressible in
$N$
and hence in
$M$.
In the latter case the subsurface of
$F$
that is contained in the solid torus
has to be injective in the solid torus. This follows
from the injectivity of 
$F$
and the fact that the intersection
$F\cap S$
has to inject into 
$M$
as well, since it is part of the pattern.
So 
$F$
consists of 
annuli in the solid torus bounded by 
$S$
(Moebius bands are ruled out by orientability of
$F$).
This means that again
$F$
has to contain at least two surface fibres of
$N$.
Take an annulus from
$F$
that is outermost in the solid torus bounded by 
$S$.
The two fibres of
$N$,
that are contained in
$F$
and lie one on each side of this annulus, can be used
to show that in this case 
$F$
has to be compressible as well. 
This proves that
the surface
$\partial N$
is incompressible in
$M$.

Since the submanifold
$N$
is disjoint from the characteristic submanifold
$\Sigma$,
each boundary torus from
$\partial N$
has to be parallel either to a canonical torus in 
$M$
or to a toral boundary component of
$M$.
Therefore 
$N$
has to be a strongly simple piece of the JSJ-decomposition
of
$M$.
If an essential annulus in
$(N, \partial N)$ 
is canonical then, by proposition 4.1 from~\cite{neumann},
it has to be a matched annulus,  making 
$N$
into a Seifert fibred space. If
$(N, \partial N)$ 
contains no canonical annulus and yet contains an 
essential one, then it is again Seifert fibred by 
proposition 3.2 in~\cite{neumann}. This implies that
$N$
is in fact an-annular.

It follows from the construction of 
$N$
that any surface from the intersection
$N\cap S_1$
will either contain the whole boundary component of
$N$
or it will contain a surface fibre of 
$N$.
The latter possibility can not occur because any component of
$S_1$
that is not contained in
$\partial N$
has to be disjoint from it and can therefore not carry a fibre of 
$N$.
The intersection
$N\cap S_2$
consists of components of
$S_2$
which are therefore surface fibres of 
$N$.
But since 
$S_2$
contains at most two parallel copies of an incompressible
$\partial$-incompressible surface the lemma follows.
\end{proof}


\vspace{-8mm}
\begin{center}
\end{center}

\begin{center}
\subsection{\normalsize \scshape TOPOLOGICAL COMPLEXITY}
\label{subsec:top}
\end{center}

There are two kinds of complexities of the surfaces 
in the canonical hierarchy we need to consider.
First there is the \textit{normal complexity},
i.e. the number of normal pieces a minimal
weight representative in the isotopy class of the surfaces consist of. 
We will deal with it in subsection~\ref{subsec:norcom}. 
Second
there is the \textit{topological complexity} of the surfaces in the
hierarchy,
that is defined in terms of their components in the following way.
To each component we assign its negative Euler characteristic and
then define the complexity 
to be the sum over all of its components. 
Topological complexity of the surfaces 
$S_1$
and
$S_2$
will, together with the Kneser-Haken finiteness 
theorem, determine the number of connected surfaces we needed to 
cut along. 
Bounding
it therefore provides a crucial step in the actual construction
of the canonical hierarchy.
Since there are
no 2-spheres, discs or projective planes in the first two steps of the
hierarchy, our topological complexity coincides with 
the Thurston complexity as defined in~\cite{thurston}.\\

\begin{lem}
\label{lem:topcomp}
Let 
$M$
be a triangulated irreducible 3-manifold 
with (possibly empty) incompressible boundary
that is different from
$\RR P^3$.
Assume further that the manifold
$M$
is atoroidal and that
it is not homeomorphic to an 
$I$-bundle. Then any
closed injective surface in
$M$
that is not boundary parallel and that has the
smallest topological complexity among all such 
surfaces is
ambient isotopic to a fundamental surface. 
\end{lem}

\begin{proof}
The general strategy of the proof is to apply a version 
of theorem~\ref{thm:sum} 
to an injective surface of minimal complexity and thus conclude that
it is fundamental. But first we need to collect some facts about
$M$
and about normal surfaces in it.

\noindent \textbf{Claim 1.} 
$M$
can not contain any of the following surfaces: 
a projective plane, 
a Moebius band, 
an injective Klein bottle. 

A boundary of a regular neighbourhood of an
injective Klein bottle is
an incompressible
torus. But all such tori in
$M$
are 
boundary parallel, making
$M$
into an 
$I$-bundle.
Since
$M$
is not homeomorphic to
$\RR P^3$,
it can not contain an embedded projective plane. The 
horizontal boundary of a regular neighbourhood of the
Moebius band is an annulus. If it is incompressible, then it
has to be 
$\partial$-parallel, which makes 
$M$
into a solid torus. This is a contradiction since
$\partial M$
is assumed to be incompressible. If our annulus is compressible,
then, the proof of the claim 
of proposition~\ref{prop:topcomp} implies that it either
cuts off a cylinder from
$M$ 
that is disjoint from the Moebius band,
or that we can find a 
3-ball that contains our Moebius band. The latter is clearly a
contradiction and the former makes 
$M$
into an 
$I$-bundle over 
$\RR P^2$.
So the claim follows.
The next sublemma is well-known.

\noindent \textbf{Sublemma.} If 
$F$
is a connected normal surface which is a sum
$F=F_1+F_2$,
then it can be expressed as another sum
$F=F_1'+F_2'$,
where both normal surfaces
$F_1'$
and
$F_2'$
are connected and no component of the intersection
$F_1'\cap F_2'$
separates both surfaces. Furthermore we can assume that
the sum 
$F=F_1'+F_2'$
is in reduced form.

The proof is by making regular alterations along certain
arcs and circles in
$F_1\cap F_2$.
The procedure, we are just about to describe, must terminate
because there is only finitely many components in
$F_1\cap F_2$.
If, say,
$F_1$
is not connected, fix one of its components and do 
all regular alterations along curves in the intersection between all 
other components and 
$F_2$.
Name the component we fixed (and didn't touch)
$F_1$
again,
and call the new normal surface 
$F_2$.
Notice that the number of components
of
$F_1\cap F_2$
is smaller than what it was before.
Now, if
$F_2$
is not connected, we can reiterate the above procedure with
$F_2$
taking the place of 
$F_1$
and vice versa.
The repetition of these steps must eventually terminate with
both
$F_1$
and
$F_2$
being connected.

Assume now that there exists a component of
$F_1\cap F_2$
that separates both surfaces into
$F_i=F_i^{(1)}\cup F_i^{(2)}$
for
$i=1,2$.
This can not be the only component of
$F_1\cap F_2$
because the surface 
$F$
is connected.
Choose notation so that the regular alteration along
this component pastes 
$F_1^{(j)}$
with
$F_2^{(j)}$
for 
$j=1,2$.
Doing regular alterations along the components of 
$F_1^{(j)}\cap F_2^{(j)}$,
for
$j=1,2$,
produces two normal surfaces, called
$F_1$
and
$F_2$
again,
with fewer components in their intersection and
with the property
$F=F_1+F_2$. 
If either of them is not connected, we repeat the
procedure for making them connected. 
Alternating between the above two processes will eventually
produce surfaces
$F_1$
and
$F_2$
that are connected, their normal sum
$F_1+F_2$
equals
$F$,
and no component of the intersection
$F_1\cap F_2$
is separating in both surfaces. 

Now take two normal surfaces
$G$
and
$H$
that satisfy all these three conditions and 
their intersection
has the smallest
number of components among all such pairs of surfaces.
We claim 
that the sum 
$F=G+H$
is in reduced form. Otherwise we can isotope the
surfaces 
$G$
and
$H$
to normal surfaces
$G'$
and
$H'$
so that the sum
$F=G'+H'$
is in reduced form. This means that the number
of components in
$G'\cap H'$
is strictly smaller than the number of pieces of
$G\cap H$.
The surfaces
$G'$
and
$H'$
are still connected, but the third condition from above must fail.
But then we can repeat the procedure described above on
$G'$
and
$H'$,
possibly reducing the number of components in
$G'\cap H'$
even further, making sure that both summands are connected
and that no component of
intersection separates both of them. This contradicts the choice
of 
$G$
and 
$H$
and hence proves the sublemma.

Now we can prove the lemma.
Suppose 
$F$
is a closed injective surface in
$M$
of the smallest topological complexity. 
It is therefore
connected and incompressible.
We isotope it into normal form so that
it is weight minimising. 
Assume now that 
$F$
is not fundamental and can therefore be expressed 
as a sum 
$F=F_1+F_2$
and that  
the conclusion of the sublemma is satisfied. 
Now we can apply the main theorem of~\cite{bart} to 
$F=F_1+F_2$,
without isotoping the summands,
because the sum is in reduced form.
This gives that the surfaces  
$F_1$
and
$F_2$
are also injective. Moreover, it follows from lemma~\ref{lem:patches}
that neither of the surfaces
$F_i$  
is a 2-sphere. 
%
Since the Euler characteristic is additive over 
normal sums and 
$M$
contains no projective planes, we now have
that 
$\chi(F)=\chi(F_1)$
and
$\chi(F_2)=0$.
Since 
$F_2$
can not be an injective Klein bottle,
by claim 1,
it can only be
a boundary parallel parallel torus. Since
$F_2$
is connected, it can contain only one copy of such
a torus. 

We now have to consider the intersection
$F_2\cap F_1$.
None of the simple closed curves from
$F_2\cap F_1$
are homotopically trivial in either of the two surfaces
because, by lemma~\ref{lem:patches}, there are no 
trivial patches.
So the space
$F_2\cap F_1$
is a 1-manifold that is homeomorphic to
a disjoint union of non-trivial parallel 
simple closed curves in the torus
$F_2$.
Let 
$X$
be the 
$(\mathrm{torus})\times I$
region between
$F_2$
and the toral boundary component of
$M$ 
that 
$F_2$
is parallel to.
Then the components of the surface
$F_1\cap X$
must be injective in 
$X$,
simply because the patches of
$F=F_1+F_2$
are injective by lemma~\ref{lem:patches}, and
$\partial X$
is incompressible in
$M$.
So, since 
$F_1\cap X$
contains no closed components, it
consists of incompressible annuli that are
disjoint from
the torus
$X\cap \partial M$.
Each such annulus must be topologically parallel
to an annulus in
$F_2$.
Let
$B$
be an outermost annular component of
$F_1\cap X$,
lying in the product region
$X$,
and let 
$A$
be the annulus in
$F_2$
that is parallel to
$B$.
There are three possible (essentially different) 
ways a normal alteration can 
act on 
$\partial A$. They are depicted by figure~\ref{fig:alter}.
\begin{figure}[!hbt]
\vspace{-2cm}
 \begin{center}
   \vspace{2cm}
    \epsfig{file=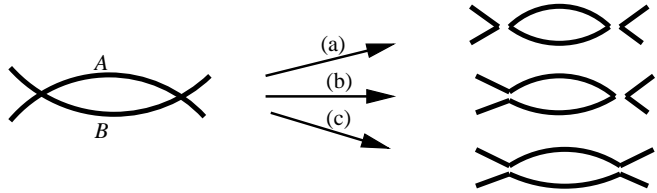}
  \caption{\small Possible normal alterations.}
  \label{fig:alter}
 \end{center}
\end{figure}

\vspace{-4mm}
They all lead to contradiction. Case (a) produces a disconnected sum.
In case (b) we can isotope the union of the patches
$A$
and
$B$
over the solid torus that they bound, to reduce the weight of
$F$.
If both 
$A$
and
$B$
had zero weight, then there would exist a normal isotopy that would
reduce the number of components in
$F_1\cap F_2$.
This contradicts the reduced form assumption. Case (c) contradicts it
as well, because the surfaces we obtain after we do the normal alterations
along
$\partial A$,
are isotopic to 
$F_1$
and
$F_2$,
but have fewer components of intersection. This proves the 
lemma.
\end{proof}

Now we are in the position to 
bound the topological complexity of surfaces
$S_1$
and
$S_2$
in the canonical hierarchy. The next corollary will
follow from the construction of the first
two surfaces in the hierarchy, from corollary 5.8 in~\cite{wang} and
lemma~\ref{lem:topcomp} and also from the bound in lemma 6.1 of~\cite{hass}.
Corollary~\ref{cor:topcomp} will be used in section~\ref{sec:proof} to
bound the number of connected surfaces in the canonical hierarchy.

\begin{cor}
\label{cor:topcomp}
Let 
$M$
be a Haken 3-manifold with (possibly empty) boundary. Let 
$T$
be its triangulation that consists of
$t$
3-simplices. 
Then the sum of the topological complexities of all
closed surfaces in
$M$,
that bound handlebodies, compression bodies and 
$I$-bundles in the complement of the canonical hierarchy after
step 2,
is 
bounded 
above by
$2^{150t}$.
\end{cor}

\begin{proof}
If
$F$
is a normal surface with respect to the triangulation
$T$,
then the complement of 
$F$
in
$M$
inherits a polyhedral structure from 
$T$.
It is obvious that there are at most 6 polyhedra
in any tetrahedron of
$T$,
lying in the complement of the normal pieces of
$F$,
that don't inherit a natural product structure of the form
$(\mathrm{triangle})\times I$
or
$(\mathrm{quadrilateral})\times I$
(see figure 4 in~\cite{mijatov}).
The complementary polyhedra with this product structure
are called \textit{parallelity regions}. We will
be interested in the components of the union of
all parallelity regions in
$M-F$
which are called
\textit{amalgams}.
They were precisely defined and studied 
in~\cite{lackenby} (sections 7 and 8). 
Notice that a vertical boundary of an amalgam (i.e. the 
part which is not contained in
$F$)
is 
a union of annuli. If they are all injective, we say that
the amalgam is \textit{maximal}. It is described
in~\cite{lackenby} how every amalgam can be extended uniquely
to the maximal one. 
A key 
property of maximal amalgams that will be crucial is the following:
if a component 
$M'$
of
$M-F$,
which is different from an 
$I$-bundle,
contains no canonical annuli, then all the maximal amalgams
in it are of the form
$(\mathrm{disc})\times I$.
It is easy to show that under these circumstances the 
3-manifold
$M'$
can be triangulated by less than
$18t$
tetrahedra (left as an exercise).

The negative Euler characteristic of
any normal surface is bounded above by four times the number
of normal discs that the surface contains. 
The topological complexity of the surfaces which bound handlebodies,
compression bodies and
$I$-bundle in the complement of the hierarchy after step 2,
increases each time we add in a closed component of 
$S_1$
and a non-annular component of
$S_2$.
The quantity from the corollary we would like to have an estimate on
is smaller than twice the sum of negative Euler 
characteristics of
$S_1$
and
$S_2$.
So all we need is the following:

\noindent {\textbf{Claim.}}   
The negative Euler characteristic of each component
of the surfaces 
$S_1$
and
$S_2$
is bounded above by
$2^{127t+19}$.

If a new component of the surface 
$S_1$
is added in a complementary piece which contains no canonical annuli,
then lemma~\ref{lem:topcomp} and the bound from 6.1 of~\cite{hass}
imply that
$4\cdot 5\cdot 18t (7\cdot 18t 2^{7\cdot 18t})<2^{127t+19}$
bounds its topological complexity (the bound is 4 times the normal
complexity of a fundamental surface in a triangulation
with
$18t$
tetrahedra and the factor
$5\cdot 18t$
bounds the number of distinct types of normal discs contained in the surface). 
The same bound works for the components of
$S_2$
because we always adjoin them to the complementary pieces with no
canonical annuli. Instead of applying lemma~\ref{lem:topcomp}
we have to use corollary 5.8
in~\cite{wang}.

If we are adding a component of
$S_1$
to a piece which contains canonical annuli, then the
process described in (1b) of the definition of the canonical hierarchy
guarantees that the negative Euler characteristic of the surface 
we are adding is bounded above by some linear function of 
$t$.
This is because the same is true for the parts of 
the boundary of the piece which 
are disjoint from the maximal amalgams contained in the piece.
This proves the claim.

We know already that the number of components of
$S_1$
is bounded above by
$2\cdot 20 t$
(number 2 accounts for the fact that sometimes we use two parallel
copies of the same component when defining 
$S_1$).
This implies that the number of non-annular components of
$S_2$
is less than 
$80t$.
This is because there is at most two non-annular components
of
$S_2$
in each complementary piece obtained by cutting along
$S_1$. 
The Euler characteristic of the surface from the corollary is 
twice the sum of the Euler characteristics of
$S_1$
and
$S_2$.
So our bound comes from the following estimate:
$2\cdot (40+80)t 2^{127t+19}<2^{150t}$
for
$t\geq 2$.
\end{proof}

\vspace{-10mm}
\begin{center}
\end{center}

\begin{center}
\subsection{\normalsize \scshape NORMAL COMPLEXITY}
\label{subsec:norcom}
\end{center}

We are now going to investigate how many normal discs are
needed to construct surfaces from the canonical hierarchy. 
For the components of
$S_1$
and
$S_2$
with negative Euler characteristic we will use similar
techniques to the ones that bounded their genus. For the
canonical surfaces we'll apply propositions~\ref{prop:cantori}
and~\ref{prop:canannuli}. 
But the non-pure surfaces from step 3 
will require some additional normal surface theory.
The estimates of the number
of normal pieces will be in terms of the number of tetrahedra in
the subdivision of the piece we are adjoining the surfaces to. 
Not subdividing and using the same triangulation to 
bound the normal complexity of several levels of the
canonical hierarchy would certainly lead to
much better bounds. But passing the information, which is
contained in the pattern, down the hierarchy
without both subdividing and keeping the pattern itself
in the 1-skeleton proved to be an insuperable task.\\

\begin{lem}
\label{lem:S_1}
Let 
$M$
be 3-manifold which either contains a closed injective surface
or has a non-trivial JSJ-decomposition. Assume also that 
$M$
has a triangulation
$T$
that consists of
$t$
tetrahedra. Then the surface 
$S_1$
from the canonical hierarchy can be isotoped into normal
form so that it contains not more than
$2^{350t^2}$
normal discs.
\end{lem}

\begin{proof}
We have already mentioned that propositions~\ref{prop:cantori}
and~\ref{prop:canannuli} together imply that all the surfaces
of the JSJ-system of 
$M$
can be put into normal form so that they consist of not more than
$2^{80t^2}$
normal discs. The rest of the components of 
$S_1$
were obtained by the recursion described in step 1 of the
canonical hierarchy. In order to bound their normal complexity
in terms of the number of tetrahedra in 
$M$,
we can use the same strategy as in lemma 2.6 in~\cite{mijatov}.
There are at most 
$20t$
non-parallel connected surfaces in
$S_1$.
The argument from lemma 2.6 in~\cite{mijatov} gives us that
$S_1$
contains not more than 
$2(2\cdot 11t 2^{11t})^{20t}$
copies of a single normal disc type from
$T$.
The factor 2 in front of the bracket is there
to account for the parallel copies of the components of
$S_1$
that we sometimes have to adjoin.
All together there are at most
$5t$
possible distinct normal disc types contained in
$S_1$.
So we have the following inequality:
$5t2(2\cdot 11t 2^{11t})^{20t}+2^{80t^2}<2^{t+4}(2^{12t+5})^{20t}+
2^{80t^2}<2^{350t^2}$.
\end{proof}

The normal complexity of the components of the surface 
$S_2$
contained in any of the pieces of
$M-\mathrm{int}(\mathcal{N}(S_1))$ 
is bounded above by
$2^{80s^2}$,
where 
$s$
is the number of 3-simplices needed to triangulate the piece.
This follows directly from proposition~\ref{prop:canannuli} for
canonical annuli
and from corollary 5.8 in~\cite{wang} for bounded surfaces
with negative Euler characteristic. We have established already
that there are at most 
$20t$
regions in the complement of
$S_1$
where we have to insert the components of
$S_2$.
Here
$t$
denotes the number of tetrahedra in the triangulation of
$M$.
So the number of normal discs in the whole of
$S_2$
is bounded above by
$20t2^{80s^2}$.

Now we need to bound the normal complexity of the surfaces
in step 3 of the canonical hierarchy. Proposition~\ref{prop:canannuli}
does that for all $P$-canonical annuli that need to be constructed
in the process. In (3b) and (3c) we might need to adjoin trivial $P$-essential
annuli of a certain kind. In the discussion that followed the proof
of proposition~\ref{prop:canannuli} we showed that such annuli are 
fundamental, provided our 3-manifold contains no non-trivial $P$-essential
annuli. All other surfaces that we have to construct have non-empty
intersection with the pattern
$P$.
Let
$\iota(F)$
be the number of points in the intersection
$F\cap P$
for any properly embedded surface
$F$
which is transverse to the pattern,
If the pattern
$P$
is contained in the 1-skeleton then the function
$\iota$
is additive over normal summation.

In order to implement step 3 we need to describe how to 
construct the non-pure spanning annuli of (3b) and also
how to find compression discs in handlebodies from (3c).
We will see that in the construction of surfaces
with non-empty intersection with the pattern, it is crucial 
that even the trivial $P$-essential annuli are well behaved
in the way that was specified in (3b) and (3c).
The strategy in this case will
be slightly different from what we usually do. Instead of proving that
our fixed surface, which minimises the intersection with the pattern,
is (almost) fundamental, we are first going to find some 
other fundamental surface 
$F$
that has the same property as our original one,
but is not necessarily isotopic to it. 
We will then calculate
$\iota(F)$.
If we express our original surface as a sum of fundamental surfaces, we 
can use the additivity of
$\iota$
to bound the number of summands in the expression. 

\begin{lem}
\label{lem:compbody}
Let 
$M$
be an irreducible bounded 3-manifold which is triangulated by 
$t$
tetrahedra. Let 
$P$
be a non-empty boundary pattern in
$\partial M$
which is contained in the 1-skeleton of the triangulation. Assume also
that every $P$-essential annulus in
$M$
is parallel to 
an annulus in the boundary of
$M$
which intersects the pattern in a disjoint union of 
homotopically non-trivial simple closed curves.

\newcounter{f}
\begin{list}{(\alph{f})}{\usecounter{f}\parsep 0mm\topsep 1mm}
\item 
Assume further that
$M$
is a
compression body with non-empty ``negative'' boundary
$\partial_-M$ and
that each toral component of
$\partial_-M$ 
is pure.
If we fix an incompressible annulus in
$M$
whose boundary circles lie in distinct components of
$\partial M$
and which has minimal intersection with the pattern
$P$,
among all annuli satisfying these conditions, 
then we can isotope it
into normal form, by an isotopy which is invariant on the pattern,
so that it consists of not more than
$2^{40t}$
normal discs. 
If such an annulus
is not separating and we take two parallel copies of it the same
bound still holds.

\item
If 
$M$
is an 
$I$-bundle over a closed surface which is 
not a torus, then an incompressible 
$\partial$-incompressible annulus 
$A$
(or its double if necessary)
in
$M$,
which has minimal intersection with the pattern 
$P$,
can be isotoped, by an isotopy that preserves the pattern,
into normal form so that it contains less than 
$2^{40t}$
normal discs.
If the base surface is a torus and
if the pattern 
$P$
consists of two non-homotopic simple closed curves, one in
each boundary component, the same bound holds.

\item
Let
$M$
be a handlebody and
let
$D$
be a 
compression disc in
$M$
which minimises 
$\iota(D)$
among all compression discs in
$M$. 
Then
we can isotope
$D$
into normal form, by an isotopy that is invariant
on the pattern, so that it contains not more than
$2^{40t}$
normal discs.
\end{list}
\end{lem}

The crucial assumption in lemma~\ref{lem:compbody}
is the one about $P$-essential annuli in
$M$.
Put differently it says that 
$M$
contains no non-trivial $P$-essential annuli and that every 
trivial one has to be parallel to an annulus in
$\partial M$
which contains the simplest possible pattern. Before we add in 
a surface with non-empty intersection with the pattern in steps
(3b) and (3c), we always make sure that this hypothesis holds for the
piece we are in by first adjoining the trivial 
$P$-essential annuli that violate this rule.
Notice also that the annulus from (a) in lemma~\ref{lem:compbody}
is automatically 
$\partial$-incompressible since its 
boundary circles lie in distinct components of
$\partial M$.
Since it is also incompressible 
it has to be a ``spanning'' annulus that was described in 
step (3b). 

It is clear that toral boundary components of genuine
compression bodies have to be contained in 
$\partial_-M$.
The fact that all such tori are pure
can be seen as follows.
Each
such torus is incompressible in the piece we are in and therefore
has to be incompressible in the ambient manifold. 
So it is either contained in
the JSJ-system or in 
the boundary of the ambient manifold. It is clearly pure in the latter
case. If it is a canonical torus, then it is there to separate
a strongly simple piece from a Seifert fibred piece. It therefore has to 
be pure again because at this stage we haven't touched the fibred pieces.
Notice also the pattern on the
$(\mathrm{torus})\times I$
pieces, which are not contained in the characteristic submanifold
$\Sigma$,
satisfies the assumption from (b) of lemma~\ref{lem:compbody}.\\

\begin{proof} We will start by proving (a).
Let 
$F$
be an incompressible annulus
in a compression body
$M$
with one of its boundary circles lying in a given component of
$\partial_-M$
and the other in
$\partial_+M$.
Assume that 
$F$
is in normal form and that it
has minimal weight
among all possible annuli that satisfy the above conditions.
Notice that we haven't stipulated anything as far as 
$\iota(F)$
is concerned. 

\noindent{\textbf{Claim 1.}} The surface $F$ is fundamental.

Assume to the contrary that 
$F=U+V$.
We can apply the sublemma from the proof of~\ref{lem:topcomp}.
So both 
$U$
and
$V$
are connected, no component of the 1-manifold
$U\cap V$
separates both of them and the sum is in reduced form.
Now we are going to go through all the potential connected surfaces
$U$
and
$V$
that satisfy the equation
$0=\chi(U)+\chi(V)$,
showing each time that we get a contradiction.

The techniques that proved lemma~\ref{lem:patches} imply that
there are no disc patches in
$F=U+V$
which are disjoint from the boundary of
$M$.
So neither of the summands is
a 2-sphere. Since
$M$
does not contain a projective plane
we can assume that
$U$
is a disc. In that case
$V$
is either a punctured torus, a punctured Moebius band or a punctured
annulus.
There are no simple closed
curves in the 1-manifold
$U\cap V$,
because 
$U$
contains no disc patches.
Therefore the boundaries of the summands must intersect. So 
the first case can not occur because it implies that both
$\partial U$
and
$\partial V$
are contained in a single component of
$\partial M$.
In the remaining two cases every outermost arc of
$U\cap V$
chops off a patch in 
$U$
which is a genuine
$\partial$-compression disc for 
$V$,
since no component of 
$U\cap V$
separates both surfaces. Let
$C$
be one such patch in
$U$
which has the smallest weight among all 
``outermost arc'' patches in
$U$.
If we 
$\partial$-compress
$V$
along 
$C$
we obtain an annulus which satisfies the same conditions
as
$F$
(it is incompressible because one of its boundary components
coincides with a boundary circle of
$F$).
Its weight is at most that of
$F$.
But since one of the parallel copies of
$C$
(that was inserted in the 
$\partial$-compression) was pasted in as 
an irregular normal alteration, there exists a weight reducing
isotopy of the annulus.
This contradicts the minimal weight 
assumption on
$F$.

Now we have to consider the case 
$\chi(U)=\chi(V)=0$.
If one of the surfaces
is closed, then the other surface is 
an annulus which is lighter than 
$F$
and has the same properties.
Both of them can not be Moebius bands because then we would
have an embedded Moebius band in 
$M$
with its boundary contained in
$\partial_-M$.
Gluing two such compression bodies 
along their ``negative'' boundaries
gives a 3-manifold which can be embedded in
$S^3$,
but which would contain an embedded Klein bottle. This argument 
also shows
that if 
$V$
is a Moebius band, then 
$\partial V$
lives in 
$\partial_+M$
and
$U$
has to be an annulus. But in this situation 
$U$
satisfies the defining conditions for 
$F$
and 
$w(U)<w(F)$.
We get the same contradiction if 
$V$
is an annulus with both of its boundary circles contained in
a single component of
$\partial M$.
The only case left is when both 
$U$
and
$V$
are annuli with boundaries lying in distinct components
of
$\partial M$.
Since
$w(U)<w(F)$
the annulus
$U$
has to be compressible in
$M$.
So it has to bound a cylinder of the form
$(\mathrm{disc})\times I$,
because its boundary circles lie in distinct components of
$\partial M$.
There are no arcs in
$U\cap V$
which are 
$\partial$-parallel in either of the annuli. This 
is because no curve from
$U\cap V$
is separating both in 
$U$
and
$V$.
Since there are no homotopically trivial simple closed
curves in
$U\cap V$,
the intersection consists of spanning arcs in both annuli.
So 
$V$
intersects the cylinder
$(\mathrm{disc})\times I$
in a collection of discs of length four.
Now we look at the disc in
$\partial M$
which is bounded by a circle from
$\partial U$.
This disc intersects 
$V$
in a collection of arcs
that decompose it into complementary regions. 
The situation for an outermost such arc is 
described by figure~\ref{fig:alter} (take 
$A$
to be a subarc in
$\partial V$
which is contained in the disc from
$\partial M$). 
By an argument, analogous to the one following figure~\ref{fig:alter},
we arrive at a contradiction which proves the claim.


Each normal disc in
$F$
intersects the 1-skeleton of the triangulation in not more than
four points. Since 
$F$
is fundamental and therefore consists of less than
$5t\cdot7t 2^{7t}$
normal discs, the number 
$\iota(F)$
is bounded above by 
$140t^22^{7t}$.
Here we are using the fact that the pattern 
$P$
lies in the 1-skeleton.

Let 
$A$
be our spanning annulus which minimises 
$\iota(A)$.
We can assume that 
$A$
is in normal form and that it
minimises the weight in its ($P$-invariant) isotopy class.

\noindent{\textbf{Claim 2.}} If
$A=X+Y$
and
$X$
is a connected surface with non-negative Euler characteristic,
then 
$\iota(X)$
is not zero.

We can assume that the sum is in reduced form.
Then by theorem~\ref{thm:sum} and lemma~\ref{lem:patches} 
we have that both 
$X$
and
$Y$
are incompressible, $P$-boundary incompressible and that there are
no disc patches. So
$X$
can not be a pure disc or a 2-sphere.
Projective planes and a Klein bottles do not embed in 
$M$
so they can not appear as summands.
If
$X$
is a torus, then, since it is incompressible, it has to be
parallel to a component of
$\partial_- M$.
By the same argument
as the one at the end of the proof of lemma~\ref{lem:topcomp} we can
conclude that the surface 
$Y$
intersects the parallelity region (between 
$X$
and the component of
$\partial_- M$)
in 
$(\mathrm{circle})\times I$.
Since there is no pattern we are free to 
isotope the circle from
$\partial A$
in this component as much as we like. This makes the summand
$X$
redundant and leads into contradiction. 
If
$X$
were a pure $P$-boundary incompressible Moebius band, then
$2X$ 
would be a $P$-essential annulus in
$M$.
Such an annulus is 
$\partial$-parallel by assumption. Since the parallelity region can not
contain a pure Moebius band, the 3-manifold
$M$
would have to be a solid torus which is a contradiction.

The last case we need to consider is when
$X$
is a $P$-essential annulus. By assumption this annulus has to 
be trivial. Since it is pure and there are no disc patches in
$A$,
the components of
$X\cap Y$
are either spanning arcs in 
$X$
or
non-trivial simple closed curves in
$X$.
The latter possibility can not occur because it would
imply that 
$\partial A=\partial X$
which is a contradiction. 
Let 
$R$
be the solid torus between 
$X$
and an annulus in
$\partial M$.
The components of
$Y\cap R$
are patches in
$A$.
They all inject into
$M$
an hence into
$R$.
None of them is homeomorphic to a Moebius 
band because 
$A$
is orientable. So the surface 
$Y\cap R$
is either a disjoint union of discs or a disjoint union 
of annuli. In the latter case each patch from
$Y\cap R$
intersect the annulus 
$X$
in at least two spanning arcs. 
Using that we can construct a pure 
$\partial$-compression disc for any patch form
$Y\cap R$
(the pattern in 
$(\partial R)-X$
consists of parallel non-trivial simple closed curves only). 
This contradicts lemma~\ref{lem:patches}. So all the 
components of the surface 
$Y\cap R$
are discs. 
Furthermore if a single patch from
$Y\cap R$
intersects 
$X$
in more than one spanning arc, then, like before, this patch
is $P$-boundary compressible. This again contradicts
lemma~\ref{lem:patches}. 
So every patch
from
$Y\cap R$
intersects each of the components of the pattern in
$(\partial R)-X$
in precisely one point. This again follows from the
$P$-boundary incompressibility of the patches of
$A$.
So we can conclude that the patches
from
$Y\cap R$,
which are all topologically discs, are in one-to-one
correspondence with the components of
$X\cap Y$
and are isotopic by an isotopy which is invariant on the pattern.
In particular each arc from
$X\cap Y$
chops of a disc in 
$Y$.
There are essentially two ways of performing
regular alterations along the arcs in
$X\cap Y$
so that the resulting normal surface does not have more components
than 
$Y$. 
Both
of these possibilities yield a surface which is isotopic to
$Y$
via an isotopy that is invariant on the pattern. Since
$A$
is connected,
$Y$
must be an embedded annulus which is isotopic to
$A$
and which satisfies the inequality
$w(Y)<w(A)$.
This contradiction proves
claim 2.

Let's express our annulus
$A$
as a sum of fundamental surfaces:
$A=k_1F_1+\cdots+k_nF_n$.
The number of summands
in this 
expression, 
that have negative
Euler characteristic,
is bounded above by the number of disc summands. All other
surfaces in the sum have non-trivial intersection with the pattern
by claim 2.
So by claim 1 there are less than 
$280t^22^{7t}$
summands in the above expression. 
The number of normal discs
in
$A$
is thus bounded above by
$280t^22^{7t}\cdot 5t7t2^{7t}<2^{18t+14}<2^{40t-1}$.
This proves (a).
Part (b) will work in the same way. Claim 1 is easier now
because the boundary of our manifold is incompressible
which means that we can apply theorem~\ref{thm:sum}. Claim 2
is true as well. The only difference is that when we are dealing 
with toral boundary, because of our assumption on the 
pattern, we can 
perform the isotopy that leads to contradiction and thus shows that
there are no toral summands in
$A$. 
All numerical bounds are the same as in (a).

In order to prove part (c) of the lemma, we need the following
analogue of claim 2.

\noindent{\textbf{Claim 3.}} Let 
$D$
be a compression disc in the handlebody 
$M$
which minimises 
$\iota(D)$
and which is in normal form and has minimal weight in its
($P$-invariant) isotopy class.
If
$D=X+Y$
and
$X$
is a connected surface with non-negative Euler characteristic,
then 
$\iota(X)$
is non-zero.

The proof of claim 3 is analogous to that of claim 2 and is
left as an exercise. 
The rest of the proof of part (c) of our lemma is identical to
what we did at the end of the proof of (a). 
\end{proof}

The last type of complementary pieces we need to construct
some normal surfaces in are
the 
$I$-bundles that arise as components of the characteristic 
submanifold
$\Sigma$.
Our starting point is a triangulated 
$I$-bundle
$M\rightarrow B$
over a possibly non-orientable bounded surface
$B$.
Let
$V$
be the vertical boundary of
$M$
with a fixed triangulation.
This
 prescribed simplicial structure on
$V$
will arise from the 
canonical hierarchy on the other side of
$V$.
The surfaces we are looking for
are vertical compression discs that will simplify 
$M$
to a 3-ball. A vertical disc in 
$M$
is a one that fibres over an arc in
$B$.
Our collection will contain 
$1-\chi(B)$
vertical compression discs and will be highly non-unique. But
we will define it in such a way, so that any two choices are related
by a homeomorphism of
$M$
which is an identity on the vertical boundary. 
In fact this homeomorphism
will come from an automorphism of the base surface
$B$,
which is fixed on the boundary
(see lemma 6.4 in~\cite{mijatov1}). 

Let 
$V_1$
be the first annulus in an ordering of the
components of
$V$
and let
$D_1,\ldots,D_n$
be a collection of disjoint vertical discs we 
want to describe
($n=1-\chi(B)$).
We are assuming that the vertical boundary
of each
$D_i$
(i.e. 
$D_i\cap V$)
intersects the annulus 
$V_1$
and that every other annulus in 
$V$
intersects precisely one compression disc. 
Let
$g$
be the genus of 
$B$,
which, in case of a 
non-orientable surface, is a maximal number of 
$\RR P^2$
summands it contains when expressed as a connected sum. 
We also stipulate that
if
$B$
is orientable
(resp. non-orientable) the first
$2g$
(resp.
$g$)
of the compression discs have their entire 
vertical boundaries contained in 
$V_1$.
The last requirement is that, even if
$B$
is not orientable, the base surface of the 
$I$-bundle 
$M-\mathrm{int}(\mathcal{N}(D_1))$
is orientable.

Now we need to make sure that the vertical boundary of the surface
$D_1\cup\ldots\cup D_n$
interacts in a prescribed way with the triangulation of 
$V$.
Choose a fibre 
$\lambda$
in each component of
$V$
which consists of the smallest number of normal arcs
with respect to the given triangulation of
$V$. 
Using these normal fibres
we can define the family
$\mathcal{F}$
of 
$2n$
fibres that will have the property 
$\mathcal{F}=V\cap(D_1\cup\ldots\cup D_n)$.
In every annulus from
$V-V_1$
the family
$\mathcal{F}$
consists of a single copy of the fibre 
$\lambda$.
In the annulus 
$V_1$
we take 
$2g+1-\chi(B)$
(resp.
$g+1-\chi(B)$)
copies of 
$\lambda$
if the base surface 
$B$
is orientable
(resp. non-orientable).
We will now give a recursive definition of the collection 
of vertical compression
discs 
(cf. subsection 6.2 in~\cite{mijatov1}). Assuming that we have already
created a subcollection
$D_1,\ldots,D_k$,
for some
$k<n$,
whose vertical boundary lies in 
$\mathcal{F}$
and which satisfies all other requirements,
we look at the 
$I$-bundle
$M_k=M-\mathrm{int}(\mathcal{N}(D_1\cup\ldots\cup D_k))$
which inherits a natural polyhedral structure 
from the original triangulation
of
$M$.
Any choice of 
the vertical disc
$D_{k+1}$
has to lie in 
$M_k$
and can be made normal with respect to this polyhedral
decomposition. It is precisely defined as follows.

\begin{prop}
\label{prop:vdisc} 
Let
$M\rightarrow B$
be a triangulated 
$I$-bundle over a (possibly non-orientable) bounded surface
$B$
and let 
$D_1,\ldots,D_k$,
where $k\in\{1,\ldots, n-1\}$,
be normal vertical compression discs as described above.
Fix two normal arcs
$e$
and
$f$
from
$\mathcal{F}\cap M_k$,
which are not contained in the union 
$D_1\cup\ldots\cup D_k$
and are supposed to be the vertical
boundary of the next disc in our collection.
Let
$D_{k+1}$
be the normal vertical compression disc whose
vertical boundary consists of
$e\cup f$
and which minimises the weight, with respect to
the polyhedral structure on
$M_k$,
among all normal vertical compression discs that carry non-trivial
elements of
$H_2(M_k,\partial M_k;\ZZ_2)$
and which are bounded by 
$e\cup f$.
Then 
$D_{k+1}$
is fundamental in
$M_k$.
\end{prop}

The proposition~\ref{prop:vdisc} is very similar 
to proposition 6.3 of~\cite{mijatov1}. 
Its proof
is simpler
because the incompressible
summands we need to deal with are never horizontal. 
Since the triangulation 
of
$M$
induces a simplicial structure on the vertical boundary
$V$,
we can use 
$\partial V$
as a boundary pattern on 
$M$
(or even on
$M_k$).
This enables us to apply our usual techniques. The proof of~\ref{prop:vdisc}
is by contradiction. We assume that our disc
$D_{k+1}$
is a sum of two normal surfaces. Then, using the familiar patch arguments,
we can show that one of the summands has to be a vertical disc
and the other one is a pure vertical annulus. 
To get a contradiction we then proceed 
in exactly the same way as in the 
proof of 6.3 in~\cite{mijatov1}. The details are left for the reader.

Making sure that, after the first compression along
$D_1$,
the base surface of the I-bundle
$M_1$
is orientable, 
is done in exactly the same way as in
subsection 6.2 of~\cite{mijatov1}.
Also 
the same bounds on the normal complexity 
apply in the setting of proposition~\ref{prop:vdisc}.
Hence our chosen family
of compression discs contains not more than
$(2\cdot 2^{11t})^{12t}<2^{80t^2}$
(for the proof see the discussion in~\cite{mijatov1},
just before lemma 6.4).
Since any automorphism of 
$B$,
that is fixed on the boundary, extends to a homeomorphism of
the 
$I$-bundle
$M$,
we can use lemma 6.4 in~\cite{mijatov1} to go between any two such families
of compression discs.

%% file: proof.tex
Now we can  
prove theorem~\ref{thm:non-fibred}. What we need to
do is to connect any two triangulations of
a given fibre-free Haken 3-manifold using Pachner moves. 
First we subdivide both triangulations so that the 
characteristic submanifold 
$\Sigma$
is triangulated by a subcomplex in
each of the subdivisions. In the strongly simple pieces of
the JSJ-decomposition
which are not contained in 
$\Sigma$, 
the subdivisions are then further simplified
using the canonical hierarchy. The gap between the triangulations
in the components of 
$\Sigma$ 
is bridged by theorem 3.1 of~\cite{mijatov1} if they are
Seifert fibred and by proposition~\ref{prop:vdisc} if they are
$I$-bundles.

Let
$M$
be a fibre-free Haken 3-manifold with a triangulation
$T$
that contains 
$t$
tetrahedra. We are now going to construct 
a ``canonical'' triangulation 
$\overline{T}$
of the complement
of the characteristic submanifold 
$\Sigma$. 
We start by subdividing 
$T$
so that the canonical hierarchy
appears as a subcomplex of this 
subdivision. 
This subcomplex induces a simplicial structure of the manifold
$M-\mathrm{int}(\Sigma)$  
which is uniquely determined by the topology of
$M$.
The new triangulation will be closely related to the canonical hierarchy
which was described in section~\ref{sec:can}. 
We know that each complementary piece of the canonical
hierarchy in  
$M-\mathrm{int}(\Sigma)$  
is a 3-ball or a solid torus with
at least one pure annulus in its boundary. 
So we can define
$\overline{T}$
to be conical in each of the 3-balls. 
Let
$K$
be the two-dimensional polyhedron
which is a union of surfaces in the canonical hierarchy.
In order to avoid confusion we should emphasise 
that 
$K$
also
contains the surfaces from
$\partial (M-\mathrm{int}(\Sigma))$.  
The complement of the singular locus of
$K$
is a disjoint union of discs and 
pure annuli which live in the boundaries
of the solid tori. The discs will be 
contained in the 2-skeleton of
$\overline{T}$
and will be triangulated as cones on their boundaries.
Since the singular locus of
$K$
is a graph which is embedded in 
the 3-manifold 
$M$,
it already has a canonical simplicial structure. 
This induces a triangulation on each of the boundaries of 
the two-dimensional faces of the polyhedron
$K$.
So by definition the triangulation
$\overline{T}$
is
uniquely determined by the canonical hierarchy in the complementary
regions which are 3-balls.

We still need to define
$\overline{T}$
in
the complementary pieces which are solid tori with pure annuli
in their boundaries.
Notice that
it follows directly from the definition of the boundary pattern
that
a pure annulus in the boundary of a complementary solid torus
can not be homotopically trivial. It therefore induces a unique
Seifert fibration of the whole piece. If
there are several pure annuli
in the boundary of a single solid torus then they must all be disjoint.
So they induce
the same Seifert fibration of the piece. Therefore
the union of all such solid tori
is a disjoint union of
Seifert fibred spaces. Furthermore there is
a simplicial structure on all boundary components of this
Seifert fibred space
which is induced by the singular locus of
$K$. 
So we can take 
$\overline{T}$
to be the simplified triangulation of our Seifert fibred
space that was defined in the proof of theorem 3.1 in~\cite{mijatov1}.

Now we need to construct 
$\overline{T}$
using Pachner moves. The starting point is the original 
triangulation
$T$
of
$M$.
Our main tool for subdividing a triangulation,
so that the subdivision contains a given normal surface in its
2-skeleton,
will be lemma 4.1 from~\cite{mijatov1}. Before we start estimating the number
of Pachner moves we need to make, we should remind ourselves that
the notation
$e^n(x)$
stands for the composition of the exponential function
$e(x)=2^x$
with itself 
$n$
times. By lemma~\ref{lem:S_1} the surface
$S_1$
consists of
$2^{350t^2}$
normal discs in the triangulation
$T$.
Making less than 
$200t2^{350t^2}<2^{360t^2}$
Pachner moves we can subdivide 
$T$
so that the subdivision
contains 
$S_1$
in its 2-skeleton.
The number of 3-simplices in the subdivision is bounded above
by 
$s=20(t+2^{350t^2})<2^{360t^2}$.
We also know (see the discussion after lemma~\ref{lem:S_1})
that the surface 
$S_2$
consists of not more than
$20t2^{80s^2}<2^{90s^2}$
normal discs in the subdivision. Applying lemma 4.1 of~\cite{mijatov1}
again, we see that 
$$200\cdot 2^{90s^2} 2^{360t^2}<2^{100s^2}<e^2(730t^2)$$
bounds the number of Pachner moves needed to construct 
the subdivision of 
$T$
which supports 
$S_1\cup S_2$
as a subcomplex.
The same expression bounds the number of tetrahedra in this 
subdivision.

In order to see how much more we need to subdivide the current
triangulation of
$M-\mathrm{int}(\Sigma)$,
if we want it to contain the polyhedron
$K$ 
in its 2-skeleton,
we have to estimate how many connected 
surfaces arise during the implementation of step 3.
The corollary~\ref{cor:topcomp} implies that the
sum of the topological complexities of all closed surfaces 
in
$M$,
which bound 
complementary pieces after the first two steps of the hierarchy, 
is bounded
above by 
$2^{150t}$.
The
$P$-canonical annuli are never 
$\partial$-parallel in the piece they appear in. The same 
is true
of the spanning annuli which feature in substep (3b).
We have already established, when we were proving that 
the canonical hierarchy has to terminate, 
that there can be at most
$9g(\partial H)$
disjoint non-parallel incompressible annuli which are not 
boundary parallel in any complementary piece 
$H$.

It is clear from the construction that
all $P$-canonical annuli that occur in a single
complementary piece
$H$
are 
disjoint and 
therefore not parallel (the annuli
we add after we've cut
$H$
for the first time
are vertical in the product structure given by 
theorem~\ref{thm:bundle}). So 
$2\cdot 9\cdot 2^{150t}$
is an upper bound on the number of such annuli in the
canonical hierarchy (the factor 2 is there because we sometimes
need to add
two parallel copies of a surface). 
The following expression
$$40t+(2-\chi(S_2))+18\cdot 2^{150t}<
 20\cdot 2^{150t}$$
bounds the number of $P$-essential annuli in the canonical
hierarchy.
We have shown (see page~\pageref{bound}) that the
first two summands control the number 
of boundary components of
$S_2$
and hence the number of trivial 
$P$-essential annuli in
$K$.
The number of $P$-canonical annuli is controlled by the 
exponential expression at the end. The inequality follows from 
the bound on
$2-\chi(S_2)$
which can be found in 
the proof of corollary~\ref{cor:topcomp}.
We can now conclude that the total number of surfaces in step 3 of the 
hierarchy is bounded above by 
$2^{160t}$.

Proposition~\ref{prop:canannuli} and lemma~\ref{lem:compbody} tell us that
if we are looking for any of the surfaces in some complementary
piece with
$r$
tetrahedra, we can construct it by using less than
$2^{80r^2}$
normal discs. By lemma 4.1 in~\cite{mijatov1} we need to make not more than
$200r2^{80r^2}$
Pachner moves to make this surface part of the 2-skeleton of the
subdivision. The same expression also bounds the number of tetrahedra
in the subdivision. It is clearly smaller than
$e^2(r)$
for 
$r$
larger than say
100.
Since the numbers we are going to apply this to are significantly bigger
than that, we can use this bound. In other words we can make 
the whole polyhedron
$K$
a subcomplex of some subdivision of the triangulation of
$M$
by making less than
$$e^{2\cdot 2^{160t}}(e^2(730t^2))<e^{2\cdot 2^{160t}}(e^3(10t))<
 e^{2^{170t}}(t)$$
moves.
Again this expression bounds the number of 3-simplices involved.
Now we have to apply theorem 5.2 
in~\cite{mijatov1} to every complementary
3-ball
region of
$K$
in order to make it conical. Theorem 3.1 
from~\cite{mijatov1} can be used to deal
with
Seifert fibred spaces which are unions of solid tori
that had pure annuli in their boundaries. Since 
$r$
(i.e. the number of tetrahedra) at this stage is so
large, the bounds in those theorems are certainly smaller than
$e^7(r)$.
The number of 3-balls and Seifert
fibred pieces we need to deal with
is bounded above by the number of 3-simplices in the subdivision.
Also the procedures simplifying the simplicial structure of the 
faces of
$K$
are linear in the number of tetrahedra of the subdivision.
The amalgamation of the edges in the singular locus of
$K$
takes linearly many steps as well. Processes very similar to 
these are described in the proof of theorem 6.5 in~\cite{mijatov1}.
So we can assume that after
$e^{2^{180t}}(t)$
Pachner moves our subdivision looks like
$\overline{T}$.

Proposition~\ref{prop:vdisc} implies that a triangulation of an
$I$-bundle over a bounded surfaces can be dealt
with in the same way theorem 6.5 in~\cite{mijatov1}
deals with 
$S^1$-bundles over bounded surfaces. 
In fact we can obtain a complete analogue of theorem 6.5 
from~\cite{mijatov1} for 
$I$-bundles over bounded surfaces. We then apply it to
the 
$I$-bundle components of
$\Sigma$.
All that is left now is to apply theorem 3.1 of~\cite{mijatov1} to the
Seifert fibred components of
$\Sigma$.
This gives the bound from theorem~\ref{thm:non-fibred}.
If our manifold is an
$I$-bundle over a closed surface
which is not Seifert fibred, we first look for 
some vertical annulus which is fundamental, and
then do the procedure described above
to its complement. Clearly the bound from theorem~\ref{thm:non-fibred}
still applies.